\documentclass[aps,pra,letterpaper,notitlepage,longbibliography,reprint,onecolumn,nofootinbib]{revtex4-2}
\makeatletter
\newcommand\orcid[1]
{
	\@AF@join{\href{https://orcid.org/#1}{ORCID:#1}}%
}%

\def\@homepage#1#2{%
	\endgroup
	\@AF@join{#1\href{https://#2}{#2}}%
}%
\makeatother

\let\oldonlinecite\onlinecite
\renewcommand{\onlinecite}[2][]{\ifblank{#1}{Ref.~[\oldonlinecite{#2}]}{#1 of Ref.~[\oldonlinecite{#2}]}}

\newcommand{\onlinecites}[2][]{Refs.~[\oldonlinecite[#1]{#2}]}
\usepackage[letterpaper,centering,hmargin=2cm,vmargin=2cm]{geometry}

\usepackage[table]{xcolor}
\usepackage{etoolbox}
\makeatletter
\patchcmd{\@hex@@Hex}{f\else}{F\else}{\typeout{Patching xcolor}}{}
\makeatother

\usepackage{environ,xstring,algpseudocode,listofitems}

\usepackage{adjustbox}
\usepackage{soul}
\usepackage{graphicx,amsmath,amsfonts,amssymb,amsthm,bm,bbm}
\usepackage{tikz,pgfplots}
\pgfplotsset{compat=newest}
\pgfplotsset{plot coordinates/math parser=false}
\usetikzlibrary{arrows,decorations.pathmorphing,backgrounds,positioning,fit,calc,patterns,shapes,decorations.markings,shadings,decorations.pathreplacing,decorations.markings,cd}
\usetikzlibrary{external}
\usepackage{ifthen}
\usepackage{multirow}
\usepackage{diagbox}
\usepackage{stmaryrd}
\usepackage{tabularx}
\usepackage{booktabs}
\newcommand{\greyhline}{\arrayrulecolor[gray]{.9}\hline\arrayrulecolor{black}}

\usepackage[pdfpagelabels,pdftex,bookmarks,breaklinks]{hyperref}
\definecolor{darkblue}{RGB}{0,0,127} 
\definecolor{darkgreen}{RGB}{0,180,0}
\definecolor{darkred}{RGB}{180,0,0}
\usepackage[capitalise]{cleveref}

\definecolor{UGentBlue}{RGB}{30,100,200}

\hypersetup{
	colorlinks,
	linkcolor=UGentBlue,
	citecolor=UGentBlue,
	filecolor=red,
	urlcolor=UGentBlue,
	pdftitle={Invertible bimodule categories and generalized Schur orthogonality},
	pdfauthor={Jacob C Bridgeman, Laurens Lootens, Frank Verstraete}
}
\let\originalleft\left
\let\originalright\right
\renewcommand{\left}{\mathopen{}\mathclose\bgroup\originalleft}
\renewcommand{\right}{\aftergroup\egroup\originalright}

\newcommand{\ZZ}[1]{\mathbb{Z}/#1\mathbb{Z}}

\DeclareMathOperator{\Tr}{Tr}

\DeclareMathOperator{\id}{id}

\newcommand{\symmetricGroup}[1]{\ensuremath{S_{#1}}}

\DeclareMathOperator{\spanop}{span}
\newcommand{\Cspan}{\spanop_{\mathbb{C}}}

\newcommand*{\one}{\mathbbm{1}} 

\newcommand{\norm}[1]{\left\lVert#1\right\rVert}

\newcommand{\restrict}[1]{\raise-.2ex\hbox{\ensuremath|}_{#1}}
\newcommand{\dualSpace}[1]{\ensuremath\hat{#1}}

\newcommand{\define}[1]{\emph{#1}}

\newcommand{\set}[2]{\left\{#1\ifthenelse{\equal{#2}{}}{}{\,\middle|\,#2}\right\}}

\newcommand{\Mat}[2][]
{
	\mathbb{M}_{#2}\ifstrempty{#1}{}{\left(#1\right)}
}

\newcommand{\dualpairing}[2]{\ensuremath{\left\langle #1, #2 \right \rangle}}

\definecolor{nicegreena}{RGB}{1,115,16}
\definecolor{nicegreenb}{RGB}{1,240,16}
\definecolor{nicegreen}{RGB}{60,183,82}

\usepackage{hhline}

\newenvironment{subalign}[1][]{\subequations\label{#1}\align}{\endalign\endsubequations}

\newboolean{default}
\newcommand{\default}{}

\newenvironment{switch}[1]
{
	\setboolean{default}{true}
	\renewcommand{\case}[2]{\ifthenelse{\equal{#1}{##1}}
		{
			\setboolean{default}{false}##2}{}
	}
	\renewcommand{\default}[1]{\ifthenelse{\boolean{default}}{##1}{}}
}{}
\tikzset{
	on each segment/.style={
		decorate,
		decoration={
			show path construction,
			moveto code={},
			lineto code={
				\path [#1]
				(\tikzinputsegmentfirst) -- (\tikzinputsegmentlast);
			},
			curveto code={
				\path [#1] (\tikzinputsegmentfirst)
				.. controls
				(\tikzinputsegmentsupporta) and (\tikzinputsegmentsupportb)
				..
				(\tikzinputsegmentlast);
			},
			closepath code={
				\path [#1]
				(\tikzinputsegmentfirst) -- (\tikzinputsegmentlast);
			},
		},
	},
	mid arrow/.style={postaction={decorate,decoration={
				markings,
				mark=at position .5 with {\arrow[#1]{stealth}}
	}}},
	mid arrow reversed/.style={postaction={decorate,decoration={
				markings,
				mark=at position .5 with {\arrowreversed[#1]{stealth}}
	}}},
	Q arrow/.style={postaction={decorate,decoration={
				markings,
				mark=at position .25 with {\arrow[#1]{stealth}}
	}}},
	TQ arrow/.style={postaction={decorate,decoration={
				markings,
				mark=at position .75 with {\arrow[#1]{stealth}}
	}}},
Q arrow reversed/.style={postaction={decorate,decoration={
			markings,
			mark=at position .25 with {\arrowreversed[#1]{stealth}}
}}},
TQ arrow reversed/.style={postaction={decorate,decoration={
			markings,
			mark=at position .75 with {\arrowreversed[#1]{stealth}}
}}},
}

\usepackage{ifplatform}

\newboolean{allowexternal}
\setboolean{allowexternal}{true}
\newboolean{onoverleaf}
\setboolean{onoverleaf}{false}

\iflinux

\def\overleafhome{/tmp}

\makeatletter
\begingroup\endlinechar=-1\relax
       \everyeof{\noexpand}%
       \edef\x{\endgroup\def\noexpand\homepath{%
		         \@@input|"kpsewhich --var-value=HOME" }}\x
\makeatother

\ifx\homepath\overleafhome

\setboolean{allowexternal}{false}
\setboolean{onoverleaf}{true}
\fi

\fi

\ifthenelse{\not\equal{\tikzexternalcheckshellescape}{}\and\boolean{allowexternal}}{
	\tikzexternalize[prefix=figures/]
	\iflinux
	\tikzset{external/system call={
			pdflatex \tikzexternalcheckshellescape -halt-on-error -interaction=batchmode -jobname "\image" "\texsource";rm -f "\image".log; rm -f "\image".dpth; rm -f "\image".md5; rm -f "\image"Notes.bib;
	}}
	\fi
	\ifwindows
	\tikzset{external/system call={
			pdflatex \tikzexternalcheckshellescape -halt-on-error -interaction=batchmode -jobname "\image" "\texsource" 
	}}
	\fi
}{
}

\newcommand{\includeTikz}[2]
{
	\tikzifexternalizing
	{
		\includeTikzrm{#1}{#2}
	}
	{
		\IfFileExists{figures/#1.pdf}{
		    \ifthenelse{\boolean{onoverleaf}}{
		    \includeTikzrm{#1}{#2}}{
			\includegraphics{figures/#1}}
		}
		{
			\includeTikzrm{#1}{#2}
		}
	}
}

\newcommand{\includeTikzrm}[2]{
	\tikzset{external/remake next}
	\tikzsetnextfilename{#1}
    #2
}

\NewEnviron{tikzarray}[2][scale=1]{
	\ensuremath{
		\begin{array}{c}
			\IfSubStr{#2}{:}
			{
				\StrBehind{#2}{:}[\filename]
				\includeTikzrm{\filename}{\begin{tikzpicture}[#1]\BODY;\end{tikzpicture}}
			}
			{
				\includeTikz{#2}{\begin{tikzpicture}[#1]\BODY;\end{tikzpicture}}
			}
	\end{array}}
}

\NewEnviron{tikzarrayrm}[2][scale=1]{
	\begin{array}{c}
		\IfSubStr{#2}{:}
		{
			\StrBehind{#2}{:}[\filename]
			\includeTikzrm{\filename}{\begin{tikzpicture}[#1]\BODY;\end{tikzpicture}}
		}
		{
			\includeTikzrm{#2}{\begin{tikzpicture}[#1]\BODY;\end{tikzpicture}}
		}
	\end{array}
}

\def\centerarc[#1](#2)(#3:#4:#5)
{ \draw[#1] ($(#2)+({#5*cos(#3)},{#5*sin(#3)})$) arc (#3:#4:#5) }

\AtBeginEnvironment{tikzcd}{\tikzexternaldisable}
\AtEndEnvironment{tikzcd}{\tikzexternalenable}
\DeclareMathOperator{\Hom}{Hom}

\DeclareMathOperator{\rk}{rk}

\newcommand{\vvec}[2][]{\operatorname{\bf Vec}^{#1}_{#2}}

\newcommand{\rrep}[1]{\operatorname{\bf Rep}\ifstrequal{#1}{}{}{\left(#1\right)}}
\newcommand{\ann}[2]{\operatorname{\bf Ann}_{#1}\left(#2\right)}
\newcommand{\Tub}[2]{\ann{#1}{#2}}

\newcommand{\bpr}[1]{\operatorname{\bf BPR}\ifstrequal{#1}{}{}{\left(#1\right)}}

\newcommand{\TY}[1]{\mathcal{TY}\ifstrequal{#1}{}{}{\left(#1\right)}}

\newcommand{\dual}[1]{\ifstrequal{#1}{}{}{\bar{#1}}}

\newcommand{\cat}[1]{\ensuremath{\mathcal{#1}}}
\newcommand{\C}{\cat{C}}
\newcommand{\M}{\cat{M}}
\newcommand{\D}{\cat{D}}

\makeatletter

\let\oldequiv\equiv

\newcommand{\@equiv}{\cong}

\renewcommand{\equiv}{\@equiv}
\newcommand{\tensorEquiv}{\ensuremath{\@equiv}}
\newcommand{\moritaEquiv}{\ensuremath{\@equiv}}

\newcommand{\iso}{\ensuremath{\simeq}}

\makeatother

\newcommand{\leftmodule}[3][]{\bimodule[#1]{#2}{#3}{}}
\newcommand{\rightmodule}[3][]{\bimodule[#1]{}{#2}{#3}}
\newcommand{\invertibleBimodule}[3][]{\bimodule[#1]{#2}{#3}{\End{#2}{#3}}}

\newcommand{\bimodule}[4][]{
	\ensuremath{
		\ifstrempty{#1}{{
				\vphantom{#3}}_{#2}#3_{#4}
			}{
				\def\txt{#3}
				\ifstrempty{#2}{}{\preto{\txt}{#2 \curvearrowright}}
				\ifstrempty{#4}{}{\appto{\txt}{\curvearrowleft #4}}
				\txt
			}		
	}
}

\definecolor{CColor}{RGB}{255,0,0}
\definecolor{MColor}{RGB}{0,50,204}
\definecolor{DColor}{RGB}{0,0,0}

\newcommand{\dimQ}[1]{d_{#1}}
\newcommand{\dimC}[1]{d_{\textcolor{CColor}{#1}}}
\newcommand{\dimM}[1]{d_{\textcolor{MColor}{#1}}}
\newcommand{\dimD}[1]{d_{\textcolor{DColor}{#1}}}
\newcommand{\FPdim}[1]{\ensuremath{\operatorname{FPdim}#1}}

\newcommand{\deltaM}[2]{\delta_{\textcolor{MColor}{#1}}^{\textcolor{MColor}{#2}}}
\newcommand{\deltaD}[2]{\delta_{\textcolor{DColor}{#1}}^{\textcolor{DColor}{#2}}}

\renewcommand{\dim}{\ensuremath{\operatorname{dim}_{\mathbb{C}}}}

\newcommand{\N}[3]{N_{#1#2}^{#3}}
\newcommand{\NCM}[3]{N_{\textcolor{CColor}{#1}\textcolor{MColor}{#2}}^{\textcolor{MColor}{#3}}}
\newcommand{\NMD}[3]{N_{\textcolor{MColor}{#1}\textcolor{DColor}{#2}}^{\textcolor{MColor}{#3}}}

\newcommand{\irr}[1]{\ensuremath{\operatorname{Irr}#1}}

\newcommand{\drinfeld}[1]{\ensuremath{\mathcal{Z}\left(#1\right)}}

\newcommand{\vecOb}{\star}

\newcommand{\CG}[7][]{
	\left(
	\begin{array}{cc|cc}
		#2 & #3 & #4 & #1\\
		#5 & #6 & \multicolumn{2}{c}{#7}
	\end{array}
	\right)
	}

\newcommand{\FLabel}[1]
{
	\begin{switch}{#1}
		\case{0}{\otimes}
		\case{1}{\triangleright}
		\case{2}{\triangleright\hspace{-.1em}\triangleleft}
		\case{3}{\triangleleft}
		\case{4}{\circ}
		\default{#1}
	\end{switch}
}

\newcommand{\FName}[1]
{
	{}^{\FLabel{#1}\!}F
}

\newcommand{\printColorList}[2]
{
	\readlist*\A{#1}
	\readlist*\B{#2}
	\def\n{\listlen\A[]}
	\value{\n}
	\foreach \i in {1,...,\n}
	{
		\textcolor{\B[\i]}{\A[\i]}
		\ifthenelse{\i=\n}{}{,}
	}
}

\makeatletter

\newcommand{\Fi}[7][]
{
	\@F[*]{#1}{#2}{#3}{#4}{#5}{#6,}{#7,}
}

\newcommand{\F}[7][]
{
	\@F[]{#1}{#2}{#3}{#4}{#5}{#6,}{#7,}
}

\newcommand{\X}[4]{
	\biggl[X_{\catOb[1]{#1}\catOb[2]{#2}}^{\catOb[1]{#3}}\biggr]_{\catOb[1]{#4}}
}

\newcommand{\catOb}[2][]
{
	\begin{switch}{#1}
		\case{0}{\textcolor{CColor}{#2}}
		\case{C}{\textcolor{CColor}{#2}}
		\case{1}{\textcolor{MColor}{#2}}
		\case{M}{\textcolor{MColor}{#2}}
		\case{2}{\textcolor{DColor}{#2}}
		\case{D}{\textcolor{DColor}{#2}}
		\default{#2}
	\end{switch}
}

\newcommand{\homSpace}[3]
{
	#1\left(#2,#3\right)
}

\newenvironment{xsmallmatrix}[1]
{\renewcommand\thickspace{\kern#1}\smallmatrix}
{\endsmallmatrix}

\newcommand{\@F}[8][] 
{
		\begin{switch}{#2}
				\case{0}{\gdef\@cA{0}\gdef\@cB{0}\gdef\@cC{0}\gdef\@cE{0}\gdef\@cF{0}\gdef\@cD{0}}
				\case{1}{\gdef\@cA{0}\gdef\@cB{0}\gdef\@cC{1}\gdef\@cE{0}\gdef\@cF{1}\gdef\@cD{1}}
				\case{2}{\gdef\@cA{0}\gdef\@cB{1}\gdef\@cC{2}\gdef\@cE{1}\gdef\@cF{1}\gdef\@cD{1}}
				\case{3}{\gdef\@cA{1}\gdef\@cB{2}\gdef\@cC{2}\gdef\@cE{1}\gdef\@cF{2}\gdef\@cD{1}}
				\case{4}{\gdef\@cA{2}\gdef\@cB{2}\gdef\@cC{2}\gdef\@cE{2}\gdef\@cF{2}\gdef\@cD{2}}
				\default{\gdef\@cA{3}\gdef\@cB{3}\gdef\@cC{3}\gdef\@cE{3}\gdef\@cF{3}\gdef\@cD{3}}
			\end{switch}
		\readlist*\@e{#7}
		\readlist*\@f{#8}
		\ifthenelse{\equal{#1}{}}
		{
			\ifthenelse{\equal{#2}{}}{\def\A{}}{\def\A{{}^{\FLabel{#2}\!}}}
		}{
				\def\A{{}_{\FLabel{#2}}}
		}
		\def\txt{\Bigl[\A F_{\catOb[\@cA]{#3}\catOb[\@cB]{#4}\catOb[\@cC]{#5}}^{\phantom{\catOb[\@cA]{#3}}\mkern-4mu\catOb[\@cD]{#6}}\Bigr]}
		\ifthenelse{\listlen\@e[]<3\AND\listlen\@f[]<3}
		{
			\preto{\txt}{
				\begin{xsmallmatrix}{.01ex}
					\catOb[\@cD]{\phantom{1}}\\
					\catOb[\@cE]{\phantom{1}}\\
					\catOb[\@cE]{\@e[1]}
				\end{xsmallmatrix}\!\!
			}
			\appto{\txt}{
				\!\!
				\begin{xsmallmatrix}{.01ex}
					\catOb[\@cD]{\phantom{1}}\\
					\catOb[\@cF]{\phantom{1}}\\
					\catOb[\@cF]{\@f[1]}
				\end{xsmallmatrix}
			}
		}
		{
			\preto{\txt}{
				\begin{xsmallmatrix}{.01ex}
					\catOb[\@cD]{\@e[3]}\\
					\catOb[\@cE]{\@e[2]}\\
					\catOb[\@cE]{\@e[1]}
				\end{xsmallmatrix}\!\!
			}
			\appto{\txt}{
				\!\!
				\begin{xsmallmatrix}{.01ex}
					\catOb[\@cD]{\@f[3]}\\
					\catOb[\@cF]{\@f[2]}\\
					\catOb[\@cF]{\@f[1]}
				\end{xsmallmatrix}
			}
		}
		\ifthenelse{\equal{#1}{*}\AND\equal{#2}{}}{\appto{\txt}{^{*}}}{}
		\txt
	}

\makeatother

\newcommand{\tubA}[3]
{
	\readlist*\@a{#1}
	\readlist*\@b{#2}
	\readlist*\@c{#3}
	\def\txt{\operatorname{\bf A}_{\catOb[1]{\@a[1]}\catOb[1]{\@a[2]}}^{\catOb[1]{\@b[1]}\catOb[1]{\@b[2]}}\!\bigg|}
	\ifthenelse{\listlen\@c[]<3}
	{
		\appto{\txt}{_{\catOb[0]{\@c[]}}}
	}
	{
		\appto{\txt}{_{(\catOb[1]{\@c[1]},\catOb[0]{\@c[2]},\catOb[1]{\@c[3]})}}
	}
	\txt
}

\newcounter{annularDiagramID}

\makeatletter

\newcommand{\tub}[3]{\@tub[annularDiagram_\arabic{annularDiagramID}]{#1}{#2}{#3}}
\newcommand{\tubrm}[3]{\@tub[:annularDiagram_\arabic{annularDiagramID}]{#1}{#2}{#3}}

\newcommand{\tubd}[3]
{
	\readlist*\@a{#1}
	\readlist*\@b{#2}
	\readlist*\@c{#3}
	\draw[MColor] (0,-1)--(0,-.2);
	\draw[MColor] (0,1)--(0,.2);
	\centerarc[black!20](0,0)(90:270:.2);\draw[black!20](0,.2)--(0,-.2);
	\centerarc[black!20](0,0)(90:270:1);
	\def\txt{
		\node[MColor,right,inner sep=1] at (0,-1) {\scriptsize\strut$\@b[1]$};
		\node[MColor,right,inner sep=1] at (0,-.2) {\scriptsize\strut$\@a[1]$};
		\node[MColor,right,inner sep=1] at (0,.2) {\scriptsize\strut$\@a[2]$};
		\node[MColor,right,inner sep=1] at (0,1) {\scriptsize\strut$\@b[2]$};
	}
	\ifthenelse{\listlen\@c[]>1}
	{
		\ifthenelse{\equal{\@c[2]}{1}}
		{
			\node[MColor,right,inner sep=1] at (0,-.5) {\scriptsize\strut$\@a[1]$};
			\node[MColor,right,inner sep=1] at (0,.5) {\scriptsize\strut$\@a[2]$};
		}{
			\draw[CColor] (0,-.6)to[out=135,in=270](-.5,0)to[out=90,in=225](0,.6);
			\node[MColor,right,inner sep=1] at (0,-.6) {\scriptsize\strut$\@c[1]$};
			\node[CColor,left,inner sep=1] at (-.45,0) {\scriptsize\strut$\@c[2]$};
			\node[MColor,right,inner sep=1] at (0,.6) {\scriptsize$\@c[3]$};
			\txt
		}
	}
	{
		\ifthenelse{\equal{\@c[1]}{1}}
		{
			\node[MColor,right,inner sep=1] at (0,-.5) {\scriptsize\strut$\@a[1]$};
			\node[MColor,right,inner sep=1] at (0,.5) {\scriptsize\strut$\@a[2]$};
		}{
			\draw[CColor] (0,-.6)to[out=135,in=270](-.5,0)to[out=90,in=225](0,.6);
			\node[CColor,left,inner sep=1] at (-.45,0) {\scriptsize\strut$\@c[1]$};
			\txt
		}
	}
}

\newcommand{\@tub}[4][]
{
	\begin{tikzarray}[scale=.55]{#1}
			\tubd{#2}{#3}{#4}
	\end{tikzarray}
	\stepcounter{annularDiagramID}
}

\makeatother

\usepackage{xifthen}    

\makeatletter

\DeclareRobustCommand\DetectUnderscore[1]{%
	\begingroup
	\protected@edef\@tempa{#1}%
	\@onelevel@sanitize\@tempa
	\expandafter\expandafter\expandafter\endgroup
	\expandafter\expandafter\expandafter\ifthenelse
	\expandafter\expandafter\expandafter{%
		\expandafter\expandafter\expandafter\isin
		\expandafter\expandafter\expandafter{%
			\expandafter\expandafter\string_%
			\expandafter}%
		\expandafter{%
			\@tempa}}{{(#1)}}{{#1}}%
}%

\newcommand{\End}[3][]{
	\ifblank{#1}{
		\ifblank{#2}
		{
			\@End{\cat{C}}{\cat{M}}
		}
		{
			\@End{#2}{#3}
		}
	}
	{
		\ifblank{#2}
		{
			\@@End{\cat{C}}{\cat{M}}
		}
		{
			\@@End{#2}{#3}
		}
	}
}
\newcommand{\@End}[2]{
	\DetectUnderscore{#1}%
	_{#2}^{*}%
}
\newcommand{\@@End}[2]{
	\ensuremath{\bf{End}_{#1}\left(#2\right)}
}
\makeatother

\newcommand{\trivalentvertex}[4]{
	\draw (0,0)--(0,1);
	\node[above,inner sep=.5] at (0,1) {\strut$#1$};
	\draw (0,0)--(-0.707,-0.707);
	\node[below,inner sep=.5] at (-0.707,-0.707) {\strut$#2$};
	\draw (0,0)--(0.707,-0.707);
	\node[below,inner sep=.5] at (0.707,-0.707) {\strut$#3$};
	\node[left] at (0,0) {\strut$#4$};
}

\newcommand{\irrepVector}[4]
{ 
	\draw[MColor] (0,-1)--(0,1) node[right,pos=0] {\scriptsize\strut$#1$} node[right,pos=1] {\scriptsize\strut$#2$};
	\draw[MColor,fill=MColor] (0,0) circle (.1);
	\node[right] at (0,0) {\scriptsize\strut$(\catOb[2]{#3},\catOb[1]{#4})$};
}
\theoremstyle{plain}
\newtheorem{theorem}{Theorem}

\newtheorem{claim}[]{Claim}

\theoremstyle{definition}
\newtheorem{definition}[]{Definition}
\newtheorem{lemma}[]{Lemma}

\newtheoremstyle{TheoremNum}
{\topsep}{\topsep}              
{\itshape}                      
{}                              
{\bfseries}                     
{.}                             
{ }                             
{\thmname{#1}\thmnote{ \bfseries #3}}
\theoremstyle{TheoremNum}

\newtheorem{thm_rep}{Theorem}

\usepackage{apptools}
\AtAppendix{
	\preto{\section}{%
		\clearpage
	}}
\AtAppendix{
	\numberwithin{theorem}{section}
}

\makeatletter
\newcommand{\@insertcomment}[5]{
	\addcontentsline{toc}{subsubsection}{\textcolor{#1}{#2:~#4}}
	\textcolor{#1}{#3#2:~#4#5}
}

\newcommand{\@inlinecomment}[3]{
	\@insertcomment{#1}{#2}{}{#3}{}
}

\newcommand{\@standardcomment}[3]{
	\@insertcomment{#1}{#2}{\begin{center}$\oldequiv\joinrel\oldequiv$}{#3}{$\oldequiv\joinrel\oldequiv$\end{center}}
}

\newcommand{\@defcommentsLong}[3]{
	\csdef{#1}##1{\@standardcomment{#2}{#3}{##1}}
	\csdef{#1inline}##1{\@inlinecomment{#2}{#3}{##1}}
}

\newcommand{\@defcomments}[1]{
	\@defcommentsLong{#1}{@#1Color}{\uppercase{#1}}
}

\definecolor{@jcbColor}{RGB}{179,0,0}
\definecolor{@llColor}{RGB}{0,150,0}
\definecolor{@fvColor}{RGB}{255, 145, 0}

\@defcomments{jcb}
\@defcomments{ll}
\@defcomments{fv}

\makeatother
\allowdisplaybreaks

\makeatletter
\renewcommand{\overset}[2]{\ensuremath{\mathop{\kern\z@\mbox{\ensuremath{#2}}}\limits^{\mbox{\scriptsize \ensuremath{#1}}}}}
\renewcommand{\underset}[2]{\ensuremath{\mathop{\kern\z@\mbox{\ensuremath{#2}}}\limits_{\mbox{\scriptsize \ensuremath{#1}}}}}
\makeatother

\begin{document}

\title{Invertible bimodule categories and generalized Schur orthogonality}
\author{Jacob C.\ Bridgeman}
\email{jcbridgeman1@gmail.com}
\homepage{jcbridgeman.bitbucket.io}
\orcid{0000-0002-5638-6681}

\author{Laurens Lootens}
\email{laurens.lootens@ugent.be}
\orcid{0000-0002-1364-4863}

\author{Frank Verstraete}
\orcid{0000-0003-0270-5592}

\affiliation{Department of Physics and Astronomy, Ghent University, Krijgslaan 281, S9, B-9000 Ghent, Belgium}

\date{\today}

\begin{abstract}
	The Schur orthogonality relations are a cornerstone in the representation theory of groups. We utilize a generalization to weak Hopf algebras to provide a new, readily verifiable condition on the skeletal data for deciding whether a given bimodule category is invertible and therefore defines a Morita equivalence. As a first application, we provide an algorithm for the construction of the full skeletal data of the invertible bimodule category associated to a given module category, which is obtained in a unitary gauge when the underlying categories are unitary. As a second application, we show that our condition for invertibility is equivalent to the notion of MPO-injectivity, thereby closing an open question concerning tensor network representations of string-net models exhibiting topological order. We discuss applications to generalized symmetries, including a generalized Wigner-Eckart theorem.
\end{abstract}

\maketitle

\vspace{-10pt}
\noindent\hfill\rule{.75\textwidth}{.5pt}\hfill
\vspace*{-10pt}
\tableofcontents
\noindent\hfill\rule{.75\textwidth}{.5pt}\hfill
\vspace{0pt}


\section{Introduction}\label{sec:introduction}

Fusion categories, and their representation theory, have a wide variety of applications, both in mathematics and physics. In particular, fusion categories are the categories of representations of (weak) Hopf algebras~\cite{Hayashi1999,Nikshych2003,Nikshych2004,Ostrik2003}, while in physics, they can be used to define a large class of 3-dimensional topological field theories~\cite{Reshetikhin1990,Reshetikhin1991,Turaev2016,Turaev2017,Turaev1992,Barrett1996} and provide a classification of 2-dimensional rational conformal field theories~\cite{Fuchs2002,Runkel2006,Froehlich2007}. It is often convenient, particularly for physical applications, to specify a fusion category in terms of its \define{skeletal data}.

Providing skeletal data is a very concrete way to specify a fusion category. It involves a finite list of numbers, which encode the associator isomorphisms of the category. This concreteness can be advantageous for various computational tasks, for example, computing the monoidal center. Additionally, it allows for verification of various properties/structures, with perhaps the most common example being unitarity of the underlying fusion category. In addition to specifying fusion categories via skeletal data, associated module/bimodule categories can be specified similarly.

A pair of fusion categories $\C$ and $\D$ are called \emph{Morita equivalent} if there exists an invertible bimodule category between them~\cite{Mueger2003}. Morita equivalent categories give equivalent Turaev-Viro invariants of 3-manifolds~\cite{Mueger2003,Turaev2010}, or, physically, Levin-Wen models in the same phase~\cite{Levin2005,Kitaev2012,Lootens2022b}. Although the output data is the same, it may be far more straightforward to compute these invariants given some inputs than other, Morita equivalent ones.

Many physical models have symmetries described by fusion categories and their (bi)modules. A class of methods based around tensor networks~\cite{Bridgeman2017b} are central to the study of these models. In \onlinecite{Sahinoglu2021}, a condition on the tensors, called \emph{MPO-injectivity}, was introduced as a necessary condition to ensure a constant bound on the ground state degeneracy. In \onlinecite{Lootens2021}, it was conjectured that this condition is equivalent to invertibility of the bimodule category describing the model.

\subsection{Central question}
In this work, we ask the following question:
`Given the data for a bimodule category, is there a simple way to check invertibility?'.

We answer affirmatively:
\begin{thm_rep}[\ref{thm:invertibility}]
	Let $\C,\,\D$ be unitary, skeletal, fusion categories, and $\bimodule{\C}{\M}{\D}$ an indecomposable, unitary, finitely semisimple, skeletal bimodule category. Then $\bimodule{\C}{\M}{\D}$ is invertible as a $(\C,\D)$ bimodule category if and only if: $\FPdim{\C} = \FPdim{\D}$ and
	\begin{align}
		\frac{1}{\rk\M}
		\sum_{\substack{\catOb[C]{a},\catOb[M]{b},\catOb[M]{d} \\
				\catOb[M]{\alpha},\catOb[M]{\beta},\catOb[M]{\mu},\catOb[M]{\nu}}}
		\frac{\dimC{a}}{\dimM{b}^2}\F[2]{a}{b}{c}{d}{\alpha,b,\mu}{\mu,d,\beta} \Fi[2]{a}{b}{c'}{d}{\alpha,b,\nu}{\nu,d,\beta}
		 & =\deltaD{c}{c^\prime},
	\end{align}
	\noindent where $\FName{2}$ denotes the bimodule associator.
\end{thm_rep}

\subsubsection*{Strategy}

There has been a great deal of prior work on characterizing Morita equivalence of classes of tensor categories, for example \onlinecites{Naidu2007,Uribe2017,Neshveyev2018}.
Although one could conceive of alternate techniques, such as explicitly computing the product of bimodule categories~\cite{Etingof2010}, each of these results is obtained by considering the \define{dual category} $\End{\C}{\M}\equiv\End[fun]{\C}{\M}$ of $\C$-module endofunctors. Given a left module category $\leftmodule{\C}{\M}$, the dual category is the unique fusion category such that $\invertibleBimodule{\C}{\M}$ is invertible. \onlinecites{Naidu2007,Uribe2017} ask: when the input category $\C$ is pointed, under what conditions the dual is also pointed. This places constraints on the skeletal data, which can be used to compute Morita equivalences in small examples~\cite{Mignard2017,Munoz2018}.

Our work follows a similar strategy. We make use of the \define{annular algebra} associated to $\leftmodule{\C}{\M}$, whose representation category is equivalent to the dual. Utilizing this representation theory, we obtain constraints on the skeletal data.

\subsection{Overview}

The remainder of this manuscript is structured as follows. In \cref{sec:preliminaries}, we provide a number of preliminaries, including definitions of the algebraic gadgets used throughout the remainder of this work. In \cref{sec:dualCat}, we discuss representations of the annular algebra, which are central to our work. In \cref{sec:invertibility}, we prove the main result of the work. In \cref{sec:orthogonality}, we discuss how Schur orthogonality is manifested for invertible bimodules.
\cref{sec:applications} contains two applications of these results: In \cref{sec:computing}, we show how the perspective taken in this work allows all skeletal data to be calculated. Some implications of orthogonality to the theory of tensor networks are discussed in \cref{sec:MPOinjectivity}. We conclude in \cref{sec:remarks}.

We provide a number of appendices to clarify various aspects of the work. A summary of the notation used in this work is provided in \cref{sec:symbols}. In \cref{app:WHA}, we review the weak Hopf algebra structure on the module annular algebra. This is essential to our main result. In \cref{app:VecG}, we discuss the special case of finite groups and their representations to connect to classical results.

\section{Preliminaries}\label{sec:preliminaries}
In this manuscript, we work exclusively over $\mathbb{C}$, all categories discussed are assumed to be finite, $\mathbb{C}$-linear and semisimple. Given a category $\cat{A}$, the space of morphisms $a\to b$ is denoted $\homSpace{\cat{A}}{a}{b}$.

\subsubsection*{Fusion category}
A \define{fusion category} $\cat{A}$ is a finitely semisimple rigid tensor category with simple unit. The category $\cat{A}$ is \define{skeletal} if each isomorphism class contains a single object. A \define{unitary} fusion category is a unitary tensor category equipped with the canonical spherical structure. From here, \emph{fusion category} will act as shorthand for unitary, skeletal fusion category.

We denote the Frobenius-Perron dimension of an object $a$ by $\dimQ{a}$, and of $\cat{A}$ by $\FPdim{\cat{A}}:=\sum_{a\in\irr{\cat{A}}} \dimQ{a}^2$, where $\irr{\cat{A}}$ is the set of simple objects. We denote the dimension of the Hom-space $\homSpace{\cat{A}}{a\otimes b}{c}$ by $\N{a}{b}{c}$, and refer to these spaces as \define{fusion spaces}.
\subsubsection*{(Bi)module category}
Given a fusion category $\C$, a (left) $\C$-module category $\leftmodule{\C}{\M}$ is a finitely semisimple, $\mathbb{C}$-linear category equipped with a functor $\triangleright:\C\times\M\to\M$, and a natural isomorphism $(-\otimes -)\triangleright- \iso -\triangleright(-\triangleright-)$, encoding associativity, obeying coherence conditions specified in \onlinecite[Def. 7.1.1]{Etingof2015}. The $\C$-module category is called \define{unitary} if it is equipped with a $C^*$-structure compatible with that of $\C$. From here, \emph{module category} will act as shorthand for unitary, skeletal module category.

Given a pair of fusion categories $\C$, $\D$, a $(\C,\D)$-bimodule category $\bimodule{\C}{\M}{\D}$ is simultaneous a left $\C$ and right $\D$ module category. Additionally, $\M$ is equipped with a natural isomorphism $(-\triangleright -)\triangleleft- \iso -\triangleright(-\triangleleft-)$ obeying coherence conditions specified in \onlinecite[Def. 7.1.7]{Etingof2015}.

The Frobenius-Perron dimension of an object $\catOb[M]{b}$ is defined as the unique positive solution to~\cite{Etingof2010}
\begin{subalign}
	\dimC{a}\dimM{b}                          & = \sum_{\catOb[M]{c}\in\irr{M}} \NCM{a}{b}{c} \dimM{c} \\
	\sum_{\catOb[M]{a}\in\irr{\M}} \dimM{a}^2 & =:\FPdim{\M}   = \FPdim{\C},
\end{subalign}
where $\NCM{a}{b}{c} := \dim \homSpace{\M}{\catOb[C]{a}\triangleright\catOb[M]{b}}{\catOb[M]{c}}$. Objects of $\C$ and $\M$ are colored red and blue respectively throughout this work, however they can also be identified from context.

We say that a bimodule category $\bimodule{\C}{\M}{\D}$ is \define{invertible}, if it is invertible as a 1-morphism in the Morita 3-category of fusion categories~\cite[Def.~4.5]{Etingof2010}. More concretely, this means there exist bimodule equivalences $\M\boxtimes_{\D}\M^{\mathrm{op}}\equiv \C,\, \M^{\mathrm{op}}\boxtimes_{\C}\M\equiv\D$, where $\M^{\mathrm{op}}$ is the $(\D,\C)$ bimodule category opposite to $\M$~\cite{Etingof2010}.

\subsubsection*{Skeletal data}
We will make extensive use of the string diagram notation for skeletal categories. We refer to \onlinecites{Kitaev2003,Turaev2017,Bonderson2007} for an introduction.

With bases fixed for all fusion spaces, the associativity isomorphisms are realized as explicit matrices. Given a bimodule $\bimodule{\C}{\M}{\D}$, we denote these by
\begin{subalign}[eqn:FMoves]
	 & \begin{tikzarray}[scale=0.5, every node/.style={scale=0.8}]{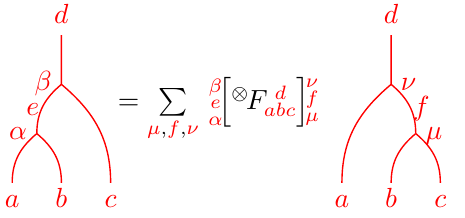}
		   \begin{scope}[local bounding box=sc1]
			\draw[CColor](0,0)to[out=90,in=220](.5,1);
			\draw[CColor](1,0)to[out=90,in=-40](.5,1);
			\draw[CColor](2,0)to[out=90,in=-40](1,2);
			\draw[CColor](.5,1)to[out=90,in=220](1,2)--(1,3);
			\node[below,inner sep=.5,CColor] at (0,0) {\strut$a$};
			\node[below,inner sep=.5,CColor] at (1,0) {\strut$b$};
			\node[below,inner sep=.5,CColor] at (2,0) {\strut$c$};
			\node[above,inner sep=.5,CColor] at (1,3) {\strut$d$};
			\node[left,CColor] at (.5,1) {\strut$\alpha$};
			\node[left,CColor] at (1,2) {\strut$\beta$};
			\node[left,CColor] at (.75,1.5) {\strut$e$};
		\end{scope}
		   \begin{scope}[local bounding box=sc2,anchor=west,shift={($(sc1.east)+(-.2,0)$)}]
			\path[] (0,-1) rectangle (7.25,1);
			\node at (0,0) {$ = \sum\limits_{\printColorList{\mu,f,\nu}{CColor,CColor,CColor}} \F[0]{a}{b}{c}{d}{\alpha,e,\beta}{\mu,f,\nu}$};
		\end{scope}
		   \begin{scope}[local bounding box=sc3,anchor=west,shift={($(sc2.east)+(-2.5,-1.5)$)}]
			\draw[CColor](0,0)to[out=90,in=220](1,2)--(1,3);
			\draw[CColor](1,0)to[out=90,in=220](1.5,1);
			\draw[CColor](2,0)to[out=90,in=-40](1.5,1);
			\draw[CColor](1.5,1)to[out=90,in=-40](1,2);
			\node[below,inner sep=.5,CColor] at (0,0) {\strut$a$};
			\node[below,inner sep=.5,CColor] at (1,0) {\strut$b$};
			\node[below,inner sep=.5,CColor] at (2,0) {\strut$c$};
			\node[above,inner sep=.5,CColor] at (1,3) {\strut$d$};
			\node[right,CColor] at (1.5,1) {\strut$\mu$};
			\node[right,CColor] at (1,2) {\strut$\nu$};
			\node[right,CColor] at (1.25,1.5) {\strut$f$};
		\end{scope}
	   \end{tikzarray},
	 & \begin{tikzarray}[scale=0.5, every node/.style={scale=0.8}]{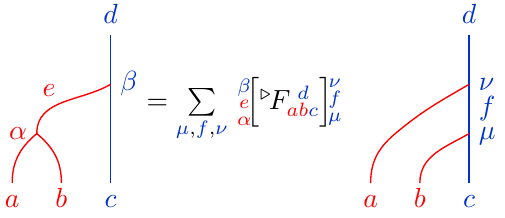}
		   \begin{scope}[local bounding box=sc1]
			\draw[MColor](2,0)--(2,3);
			\draw[CColor](1,0)to[out=90,in=-40](.5,1);
			\draw[CColor](0,0)to[out=90,in=220](.5,1)to[out=90,in=210](2,2);
			\node[below,inner sep=.5,CColor] at (0,0) {\strut$a$};
			\node[below,inner sep=.5,CColor] at (1,0) {\strut$b$};
			\node[below,inner sep=.5,MColor] at (2,0) {\strut$c$};
			\node[above,inner sep=.5,MColor] at (2,3) {\strut$d$};
			\node[left,CColor] at (.5,1) {\strut$\alpha$};
			\node[above,inner sep=1,CColor] at (.75,1.5) {\strut$e$};
			\node[right,MColor] at (2,2) {\strut$\beta$};
		\end{scope}
		   \begin{scope}[local bounding box=sc2,anchor=west,shift={($(sc1.east)+(-.2,0)$)}]
			\path[] (0,-1) rectangle (7.25,1);
			\node at (0,0) {$=\sum\limits_{\printColorList{\mu,f,\nu}{MColor,MColor,MColor}} \F[1]{a}{b}{c}{d}{\alpha,e,\beta}{\mu,f,\nu}$};
		\end{scope}
		   \begin{scope}[local bounding box=sc3,anchor=west,shift={($(sc2.east)+(-2.5,-1.5)$)}]
			\draw[MColor](2,0)--(2,3);
			\draw[CColor](1,0)to[out=90,in=210](2,1);
			\draw[CColor](0,0)to[out=90,in=220](.5,1)to[out=40,in=210](2,2);
			\node[below,inner sep=.5,CColor] at (0,0) {\strut$a$};
			\node[below,inner sep=.5,CColor] at (1,0) {\strut$b$};
			\node[below,inner sep=.5,MColor] at (2,0) {\strut$c$};
			\node[above,inner sep=.5,MColor] at (2,3) {\strut$d$};
			\node[right,MColor] at (2,1) {\strut$\mu$};
			\node[right,MColor] at (2,1.5) {\strut$f$};
			\node[right,MColor] at (2,2) {\strut$\nu$};
		\end{scope}
	   \end{tikzarray}
	,
	 & \begin{tikzarray}[scale=0.5, every node/.style={scale=0.8}]{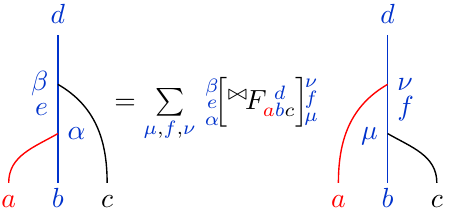}
		   \begin{scope}[local bounding box=sc1]
			\draw[MColor](1,0)--(1,3);
			\draw[CColor](0,0)to[out=90,in=210](1,1);
			\draw[DColor](2,0)to[out=90,in=-30](1,2);
			\node[below,inner sep=.5,CColor] at (0,0) {\strut$a$};
			\node[below,inner sep=.5,MColor] at (1,0) {\strut$b$};
			\node[below,inner sep=.5,DColor] at (2,0) {\strut$c$};
			\node[above,inner sep=.5,MColor] at (1,3) {\strut$d$};
			\node[right,MColor] at (1,1) {\strut$\alpha$};
			\node[left,MColor] at (1,1.5) {\strut$e$};
			\node[left,MColor] at (1,2) {\strut$\beta$};
		\end{scope}
		   \begin{scope}[local bounding box=sc2,anchor=west,shift={($(sc1.east)+(-.2,0)$)}]
			\path[] (0,-1) rectangle (7.25,1);
			\node at (0,0) {$=\sum\limits_{\printColorList{\mu,f,\nu}{MColor,MColor,MColor}} \F[2]{a}{b}{c}{d}{\alpha,e,\beta}{\mu,f,\nu}$};
		\end{scope}
		   \begin{scope}[local bounding box=sc3,anchor=west,shift={($(sc2.east)+(-2.5,-1.5)$)}]
			\draw[MColor](1,0)--(1,3);
			\draw[CColor](0,0)to[out=90,in=210](1,2);
			\draw[DColor](2,0)to[out=90,in=-30](1,1);
			\node[below,inner sep=.5,CColor] at (0,0) {\strut$a$};
			\node[below,inner sep=.5,MColor] at (1,0) {\strut$b$};
			\node[below,inner sep=.5,DColor] at (2,0) {\strut$c$};
			\node[above,inner sep=.5,MColor] at (1,3) {\strut$d$};
			\node[left,MColor] at (1,1) {\strut$\mu$};
			\node[right,MColor] at (1,1.5) {\strut$f$};
			\node[right,MColor] at (1,2) {\strut$\nu$};
		\end{scope}
	   \end{tikzarray},
	\\
	 & \begin{tikzarray}[scale=0.5, every node/.style={scale=0.8}]{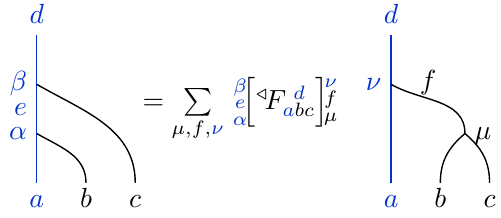}
		   \begin{scope}[local bounding box=sc1]
			\draw[MColor](0,0)--(0,3);
			\draw[DColor](1,0)to[out=90,in=-30](0,1);
			\draw[DColor](2,0)to[out=90,in=-30](0,2);
			\node[below,inner sep=.5,MColor] at (0,0) {\strut$a$};
			\node[below,inner sep=.5,DColor] at (1,0) {\strut$b$};
			\node[below,inner sep=.5,DColor] at (2,0) {\strut$c$};
			\node[above,inner sep=.5,MColor] at (0,3) {\strut$d$};
			\node[left,MColor] at (0,1) {\strut$\alpha$};
			\node[left,MColor] at (0,1.5) {\strut$e$};
			\node[left,MColor] at (0,2) {\strut$\beta$};
		\end{scope}
		   \begin{scope}[local bounding box=sc2,anchor=west,shift={($(sc1.east)+(-.2,0)$)}]
			\path[] (0,-1) rectangle (7.25,1);
			\node at (0,0) {$=\sum\limits_{\printColorList{\mu,f,\nu}{DColor,DColor,MColor}} \F[3]{a}{b}{c}{d}{\alpha,e,\beta}{\mu,f,\nu}$};
		\end{scope}
		   \begin{scope}[local bounding box=sc3,anchor=west,shift={($(sc2.east)+(-2,-1.5)$)}]
			\draw[MColor](0,0)--(0,3);
			\draw[DColor](1,0)to[out=90,in=220](1.5,1);
			\draw[DColor](2,0)to[out=90,in=-40](1.5,1)to[out=90,in=-30](0,2);
			\node[below,inner sep=.5,MColor] at (0,0) {\strut$a$};
			\node[below,inner sep=.5,DColor] at (1,0) {\strut$b$};
			\node[below,inner sep=.5,DColor] at (2,0) {\strut$c$};
			\node[above,inner sep=.5,MColor] at (0,3) {\strut$d$};
			\node[right,DColor] at (1.5,1) {\strut$\mu$};
			\node[above,DColor] at (.75,1.5) {\strut$f$};
			\node[left,MColor] at (0,2) {\strut$\nu$};
		\end{scope}
	   \end{tikzarray},
	 & \begin{tikzarray}[scale=0.5, every node/.style={scale=0.8}]{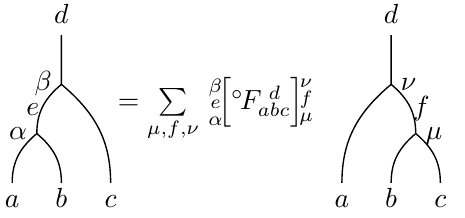}
		   \begin{scope}[local bounding box=sc1]
			\draw[DColor](0,0)to[out=90,in=220](.5,1);
			\draw[DColor](1,0)to[out=90,in=-40](.5,1);
			\draw[DColor](2,0)to[out=90,in=-40](1,2);
			\draw[DColor](.5,1)to[out=90,in=220](1,2)--(1,3);
			\node[below,inner sep=.5,DColor] at (0,0) {\strut$a$};
			\node[below,inner sep=.5,DColor] at (1,0) {\strut$b$};
			\node[below,inner sep=.5,DColor] at (2,0) {\strut$c$};
			\node[above,inner sep=.5,DColor] at (1,3) {\strut$d$};
			\node[left,DColor] at (.5,1) {\strut$\alpha$};
			\node[left,DColor] at (1,2) {\strut$\beta$};
			\node[left,DColor] at (.75,1.5) {\strut$e$};
		\end{scope}
		   \begin{scope}[local bounding box=sc2,anchor=west,shift={($(sc1.east)+(-.2,0)$)}]
			\path[] (0,-1) rectangle (7.25,1);
			\node at (0,0) {$ = \sum\limits_{\printColorList{\mu,f,\nu}{DColor,DColor,DColor}} \F[4]{a}{b}{c}{d}{\alpha,e,\beta}{\mu,f,\nu}$};
		\end{scope}
		   \begin{scope}[local bounding box=sc3,anchor=west,shift={($(sc2.east)+(-2.5,-1.5)$)}]
			\draw[DColor](0,0)to[out=90,in=220](1,2)--(1,3);
			\draw[DColor](1,0)to[out=90,in=220](1.5,1);
			\draw[DColor](2,0)to[out=90,in=-40](1.5,1);
			\draw[DColor](1.5,1)to[out=90,in=-40](1,2);
			\node[below,inner sep=.5,DColor] at (0,0) {\strut$a$};
			\node[below,inner sep=.5,DColor] at (1,0) {\strut$b$};
			\node[below,inner sep=.5,DColor] at (2,0) {\strut$c$};
			\node[above,inner sep=.5,DColor] at (1,3) {\strut$d$};
			\node[right,DColor] at (1.5,1) {\strut$\mu$};
			\node[right,DColor] at (1,2) {\strut$\nu$};
			\node[right,DColor] at (1.25,1.5) {\strut$f$};
		\end{scope}
	   \end{tikzarray},
\end{subalign}
where Greek letters denote basis vectors in the appropriate Hom-space, and we refer to the associator matrices as $F$-symbols. A bimodule category is specified by supplying the full set of \define{skeletal data}
\begin{align}
	\left(\FName{0},\FName{1},\FName{2},\FName{3},\FName{4}\right).
\end{align}

Given an $F$-symbol, we denote the inverse by lowering the label, for example
\begin{align}
	\sum\limits_{\printColorList{\mu,f,\nu}{MColor,MColor,MColor}}\F[2]{a}{b}{c}{d}{\alpha,e,\beta}{\mu,f,\nu} \Fi[2]{a}{b}{c}{d}{\alpha^\prime,e^\prime,\beta^\prime}{\mu,f,\nu}
	 & =
	\delta_{\catOb[M]{\alpha}}^{\catOb[M]{\alpha^\prime}}\delta_{\catOb[M]{e}}^{\catOb[M]{e^\prime}}\delta_{\catOb[M]{\beta}}^{\catOb[M]{\beta^\prime}}.
\end{align}

Finally, we note that there is a freedom, referred to as \define{gauge freedom}, in specifying $F$-symbols arising from changing basis on the Hom-spaces. We will always choose a basis in which the matrices are unitary (i.e. the inverse is the conjugate transpose), and the other gauge conditions outlined in \onlinecite{Barter2022} also hold.

\subsubsection*{Module annular algebra}
The algebraic gadget central to our work is the \define{module annular algebra}~\cite{Kitaev2012,Lan2014,Bridgeman2019,Hoek2019,Bridgeman2020,Barter2022}.
\begin{definition}[Module annular algebra]
	Given a (unitary) module category $\leftmodule{\C}{\M}$, the associated annular algebra $\ann{\C}{\M}$ is the algebra with basis consisting of diagrams (up to planar isotopy)
	\begin{align}
		\set{
			\begin{tikzarray}[scale=.65, every node/.style={scale=0.8}]{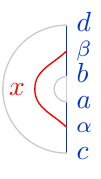}
				\draw[MColor] (0,-1)--(0,-.2);
				\draw[MColor] (0,1)--(0,.2);
				\centerarc[black!20](0,0)(90:270:.2);\draw[black!20](0,.2)--(0,-.2);
				\centerarc[black!20](0,0)(90:270:1);
				\draw[CColor] (0,-.6)to[out=135,in=270](-.5,0)to[out=90,in=225](0,.6);
				\node[MColor,right] at (0,-1) {\strut$c$};
				\node[MColor,right] at (0,-.2) {\strut$a$};
				\node[MColor,right] at (0,.2) {\strut$b$};
				\node[MColor,right] at (0,1) {\strut$d$};
				\node[MColor,right] at (0,-.6) {\footnotesize\strut$\alpha$};
				\node[MColor,right] at (0,.6) {\footnotesize$\beta$};
				\node[CColor,left] at (-.5,0) {\strut$x$};
			\end{tikzarray}
		}{\catOb[C]{x}\in\irr{\C},\catOb[M]{a},\catOb[M]{b},\catOb[M]{c},\catOb[M]{d}\in\irr{\M},1\leq\textcolor{MColor}{\alpha}\leq \NCM{x}{a}{c},1\leq\textcolor{MColor}{\beta}\leq \NCM{x}{n}{d}},\label{eqn:pictureBasis}
	\end{align}
	where $\irr{}$ denotes the set of simple objects, and $\NCM{a}{b}{c}$ is the dimension of $\homSpace{\M}{\catOb[C]{a}\triangleright\catOb[M]{b}}{\catOb[M]{c}}$. The product is given by concentric stacking, with the local relations \cref{eqn:FMoves} being used to reduce the result to the \define{picture basis} defined in \cref{eqn:pictureBasis}.

\end{definition}

\subsubsection*{Weak Hopf algebras}
In addition to being associative algebras, annular algebras can be equipped with the additional structure of a pure $C^*$-weak Hopf algebra (WHA). Since this structure will be an important part of this work, we present the action of the WHA maps in \cref{app:WHA}, and here review some important properties we will make use of.
We refer to \onlinecite{Boehm1999} for a more complete introduction to WHAs.

A weak bialgebra with underlying (finite dimensional) vector space $A$ is simultaneously an associative algebra $(A,\mu,\eta)$, and a coassociative coalgebra $(A,\Delta,\varepsilon)$. The bialgebra is only `weak' since the compatibility with the unit and counit are relaxed, see \onlinecites{Boehm1999,Nikshych2004}. In addition, a $C^*$-WHA is equipped with a $*$-structure, and an anti-isomorphism $S$ called the antipode. We refer to \onlinecite{Boehm1999} for the compatibility relations on these maps. Where it does not cause confusion, we will refer to $A$ as a WHA. Finally, by an $A$-module, we mean a module over the algebra $(A,\mu,\eta)$ for which the unit acts as the identity.

We will need several results concerning $C^*$-WHAs and their representations. We refer to \onlinecites{Boehm1999,Boehm2000,Nikshych2004} for more details. In the following, we label the trivial representation by $1$. In the remainder of this manuscript, we assume all WHAs are \emph{pure} (also called \emph{connected}~\cite{Nikshych2004}), meaning the trivial representation is indecomposable. Additionally, we assume they are finite and semisimple.

The categories of representations of semisimple weak Hopf algebras are rigid tensor categories. In the finite, pure case, these are fusion categories~\cite{Nikshych2004}. Given an irreducible representation $a$, we will denote the dual representation by $\bar{a}$, since it corresponds to the dual object in the associated fusion category.

In the following, $A$ is a finite, semisimple $C^*$-WHA.
\begin{claim}
	Let $\chi_a\in \dualSpace{A}$  be the irreducible characters of $A$. Then $\chi^*_a=\chi_{\bar{a}}$.
	\begin{proof}
		\begin{align}
			\chi_a^*(x) & :=\overline{\chi_a\left(S(x)^*\right)}=\chi_a\left(S(x)\right)=\chi_{\dual{a}}\left(x\right),
		\end{align}
		where the second equality follows from $a$ labeling a $*$-representation, so $\langle v,S(x)^* \cdot w\rangle_{V_a}=\langle S(x) \cdot v, w\rangle_{V_a}=\overline{\langle w,S(x) \cdot v\rangle}_{V_a}$. The final equality arises from \onlinecite[Section 2.2]{Boehm2000}.
	\end{proof}
\end{claim}

As a straightforward consequence of \onlinecite[Lemma 4.8]{Boehm1999}, we have
\begin{claim}\label{result:chiLambda}
	Denote the Haar integral of $A$ by $\Lambda$. Let $a$ label a simple, finite dimensional $A$-module, then
	\begin{align}
		\chi_{a}(\Lambda) & = \delta_{a,1}.
	\end{align}
\end{claim}

The coproduct gives a natural action of a WHA on the tensor product of $A$-modules. To ensure the action is non-degenerate (i.e. the unit acts as the identity), the tensor product of modules is defined as~\cite{Boehm2000}
\begin{align}
	V\boxtimes W := & \set{x\in V\otimes W}{\Delta(1)x=x}.
\end{align}

Since the category of representations of a WHA is a fusion category, we can decompose the tensor product into irreducible representations
\begin{align}
	V_{a}\boxtimes V_{b} \iso \oplus_{c\in\irr{\rrep{A}}} \N{a}{b}{c} V_c,
\end{align}
where the $\N{a}{b}{c}$ are the multiplicities. The trivial representation (denoted $1$) is the tensor unit, and each irrep $a$ has a dual $\dual{a}$ such that $\N{a}{b}{1} = \N{b}{a}{1} = \delta_{b}^{\dual{a}}$.

\section{Representations and the dual category}\label{sec:dualCat}

Given any bimodule category $\bimodule{\C}{\M}{\D}$, we can construct representations of $\ann{\C}{\M}$ as follows: For each object $\catOb[D]{b}\in\D$, we can construct a (graded) vector space
\begin{align}
	V_{\catOb[D]{b}} & := \bigoplus_{\catOb[M]{a},\catOb[M]{c}\in \irr{\M}} \homSpace{\M}{\catOb[M]{a}\triangleleft \catOb[D]{b}}{\catOb[M]{c}}.
\end{align}
It will be convenient to work with two graded bases for $V_{\catOb[D]{b}}$, related by a rescaling
\begin{align}
	\begin{tikzarray}[scale=.4]{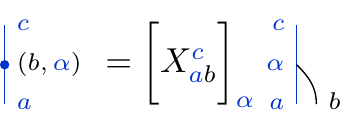}
		\begin{scope}[local bounding box=sc1]
			\irrepVector{a}{c}{b}{\alpha}
		\end{scope}
		\begin{scope}[local bounding box=sc2,anchor=west,shift={($(sc1.east)+(0,0)$)}]
			\node at (0,0) {$= \X{a}{b}{c}{\alpha}$};
		\end{scope}
		\begin{scope}[local bounding box=sc3,anchor=west,shift={($(sc2.east)+(.75,0)$)}]
			\draw[MColor] (0,-1)--(0,1) node[left,pos=0] {\scriptsize\strut$a$} node[left,pos=1] {\scriptsize\strut$c$};
			\draw[DColor] (.5,-1) to[out=90,in=-45] (0,0);\node[DColor,right] at (.5,-1) {\scriptsize\strut$b$};
			\node[left] at (0,0) {\scriptsize\strut$\catOb[1]{\alpha}$};
		\end{scope}
	\end{tikzarray}.\label{eqn:relateBases}
\end{align}
We make this distinction because the `string diagram basis' on the right side of \cref{eqn:relateBases} is part of the data of $\bimodule{\C}{\M}{\D}$, while for the purposes of representations, it is convenient to use an orthonormal basis referred to as the `dot basis'.
From this identification, we can define the action of $\ann{\C}{\M}$ on $V_{\catOb[D]{c}}$:
\begin{align}
	\begin{tikzarray}[]{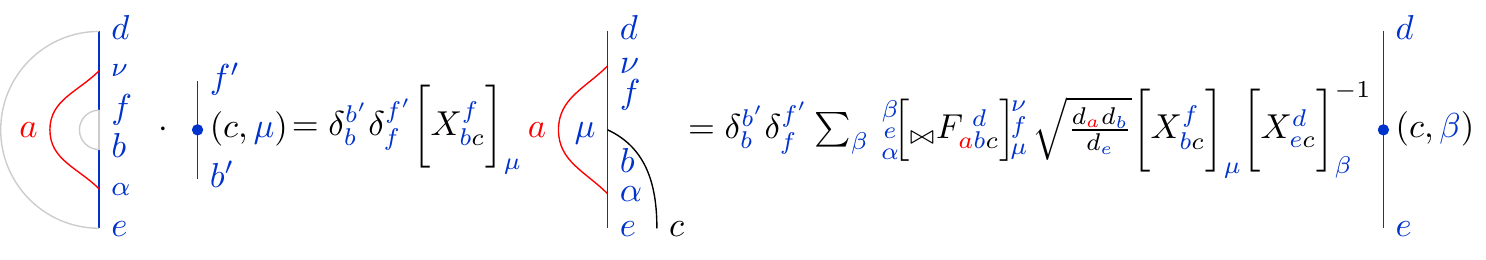}
		\pgfmathsetmacro{\s}{1/sqrt(2)};
		\begin{scope}[local bounding box=sc1]
			\draw[MColor] (0,-1)--(0,-.2);
			\draw[MColor] (0,1)--(0,.2);
			\centerarc[black!20](0,0)(90:270:.2);\draw[black!20](0,.2)--(0,-.2);
			\centerarc[black!20](0,0)(90:270:1);
			\draw[CColor] (0,-.6)to[out=135,in=270](-.5,0)to[out=90,in=225](0,.6);
			\node[MColor,right] at (0,-1) {\strut$e$};
			\node[MColor,right] at (0,-.2) {\strut$b$};
			\node[MColor,right] at (0,.2) {\strut$f$};
			\node[MColor,right] at (0,1) {\strut$d$};
			\node[MColor,right] at (0,-.6) {\footnotesize\strut$\alpha$};
			\node[MColor,right] at (0,.6) {\footnotesize$\nu$};
			\node[CColor,left] at (-.5,0) {\strut$a$};
			\node at (.65,0) {$\cdot$};
			\draw[MColor] (1,-.5)--(1,.5) node[right,pos=0] {\strut$b^\prime$} node[right,pos=1] {\strut$f^\prime$};
			\draw[MColor,fill=MColor] (1,0) circle (.05);
			\node[right] at (1,0) {\strut$(\catOb[2]{c},\catOb[1]{\mu})$};
		\end{scope}
		\begin{scope}[local bounding box=sc2,anchor=west,shift={($(sc1.east)+(-.2,0)$)}]
			\node at (0,0) {$= \deltaM{b}{b^\prime}\deltaM{f}{f^\prime}\X{b}{c}{f}{\mu}$};
		\end{scope}
		\begin{scope}[local bounding box=sc3,anchor=west,shift={($(sc2.east)+(.75,0)$)}]
			\draw[CColor] (0,-.65)to[out=135,in=270] (-.5,0) to[out=90,in=225] (0,.65);
			\draw[MColor] (0,-1)--(0,1) node[right,pos=0] {\strut$e$} node[right,pos=1] {\strut$d$};
			\draw[DColor] (.5,-1) to[out=90,in=-25] (0,0);\node[DColor,right] at (.5,-1) {\strut$c$};
			\node[left] at (0,0) {\strut$\catOb[1]{\mu}$};
			\node[right] at (0,-.65) {\strut$\catOb[1]{\alpha}$};
			\node[right] at (0,.65) {\strut$\catOb[1]{\nu}$};
			\node[right] at (0,-.35) {\strut$\catOb[1]{b}$};
			\node[right] at (0,.35) {\strut$\catOb[1]{f}$};
			\node[left] at (-.5,0) {\strut$\catOb[0]{a}$};
		\end{scope}
		\begin{scope}[local bounding box=sc4,anchor=west,shift={($(sc3.east)+(-.2,0)$)}]
			\node at (0,0) {$= \deltaM{b}{b^\prime}\deltaM{f}{f^\prime}\sum_{\catOb[1]{\beta}}\Fi[2]{a}{b}{c}{d}{\alpha,e,\beta}{\mu,f,\nu}\sqrt{\frac{\dimC{a}\dimM{b}}{\dimM{e}}}\X{b}{c}{f}{\mu}\X{e}{c}{d}{\beta}^{-1}$};
		\end{scope}
		\begin{scope}[local bounding box=sc5,anchor=west,shift={($(sc4.east)+(0,0)$)}]
			\draw[MColor] (0,-1)--(0,1) node[right,pos=0] {\strut$e$} node[right,pos=1] {\strut$d$};
			\draw[MColor,fill=MColor] (0,0) circle (.05);
			\node[right] at (0,0) {\strut$(\catOb[2]{c},\catOb[1]{\beta})$};
		\end{scope}
	\end{tikzarray}.\label{eqn:action}
\end{align}
In the representation $V$, we will denote the matrix elements of a diagram $T$ by $\rho_{V}(T)_{(\printColorList{e',\beta,d'}{MColor,MColor,MColor}),(\printColorList{b',\mu,f'}{MColor,MColor,MColor})}$.

If the category $\D$ is chosen to be (equivalent to) the category of representations of $\ann{\C}{\M}$, we denote it by $\End{\C}{\M}$, and refer to it as the \define{dual}~\cite{Ostrik2003,Etingof2015}. Choosing $\D\tensorEquiv \End{\C}{\M}$ is the unique way to extend $\leftmodule{\C}{\M}$ to an \emph{invertible} bimodule category. In that case, irreducible representations are labeled by simple objects, and morphisms in $\End{\C}{\M}$ are given by $\ann{\C}{\M}$-module maps.

The tensor product $\boxtimes$ defined on representations corresponds to vertical stacking of vectors, with $\Delta(1)x = x$ forcing matching of the common label, for example
\begin{align}
	\begin{tikzarray}[scale=.3]{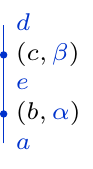}
		\begin{scope}[shift={(0,-1)}]
			\irrepVector{a}{e}{b}{\alpha}
		\end{scope}
		\begin{scope}[shift={(0,1)}]
			\irrepVector{}{d}{c}{\beta}
		\end{scope}
	\end{tikzarray}
	\in V_{\catOb[D]{b}}\boxtimes V_{\catOb[D]{c}}.
\end{align}
When $\D\tensorEquiv \End{\C}{\M}$, this tensor product will decompose into simple modules, labeled by simple objects in $\End{\C}{\M}$.
The structure constants of $V_{\catOb[D]{b}}\boxtimes V_{\catOb[D]{c}}$ can be computed using the data of $\D$ using \cref{eqn:relateBases}, with the associators $\FName{3}$ intertwining between the product representation and its decomposition in terms of irreps. Since $\bimodule{\C}{\M}{\D}$ was a bimodule category, the product $\boxtimes$ is automatically compatible with the product in $\D$. We therefore say that $\FName{3}$ isomorphisms witness that this is a tensor equivalence. In \cref{sec:computing}, we utilize this to compute $\FName{3}$.

Our main results will be based around generalizing/applying Schur's orthogonality theorems to these annular algebras. The first orthogonality theorem involves irreducible characters.
The character of an annular diagram in the module labeled by $\catOb[D]{c}\in \D$, obtained by tracing over \cref{eqn:action}, is
\begin{align}
	\chi_{\catOb[D]{c}}\left(\tub{b,f}{e,d}{\alpha,a,\nu}\right) & =\deltaM{b}{e} \deltaM{d}{f} \sqrt{\dimC{a}} \sum_{\zeta}^{\NMD{b}{c}{d}} \Fi[2]{a}{b}{c}{d}{\alpha,b,\zeta}{\zeta,d,\nu},\label{eqn:characterBasisVector}
\end{align}
which is an irreducible character exactly when $\D\tensorEquiv \End{\C}{\M}$.


\section{Invertibility of bimodules}\label{sec:invertibility}

In this section, we prove our main result. Taking inspiration from Schur's orthogonality theorems, we proceed by introducing an inner product on the space of characters. After arguing that irreducible characters are orthonormal with respect to this inner product, we show how this allows us to prove, or rule out, invertibility in the case the WHA is an annular algebra.

We begin with the inner product:
\begin{lemma}[Orthogonality of characters]\label{result:SchurOrthogonality}
	Let $A$ be a finitely semisimple, pure $C^*$-weak Hopf algebra with Haar integral $\Lambda$. Let $a$ and $b$ label simple, finite dimensional $A$-modules, then
	\begin{align}
		(\chi_a,\chi_b) & :=\dualpairing{\chi_a^*\chi_b}{\Lambda}=\delta_{a}^{b}.
	\end{align}
	\begin{proof}
		\begin{align}
			\chi_a^*\chi_b & = \chi_{\dual{a}}\chi_b = \sum_{c} \N{\dual{a}}{b}{c} \chi_c,
		\end{align}
		where $\N{a}{b}{c}$ is the multiplicity of the simple module $V_c$ in $V_a\boxtimes V_b$. We therefore have
		\begin{subalign}
			(\chi_a,\chi_b) & = \sum_{c} \N{\dual{a}}{b}{c} \chi_c (\Lambda)                 \\
			                & = \sum_{c} \N{\dual{a}}{b}{c} \delta_{c}^{1} = \delta_{a}^{b},
		\end{subalign}
		where the second equality follows from \cref{result:chiLambda}.
	\end{proof}
\end{lemma}

In the case of an annular algebra, \cref{result:SchurOrthogonality} allows us to prove, or rule out, invertibility of the bimodule.
\begin{theorem}[Invertibility]\label{thm:invertibility}
	Let $\C,\,\D$ be unitary, fusion categories, and $\bimodule{\C}{\M}{\D}$ an indecomposable, unitary, finitely semisimple, skeletal bimodule category. Then $\M$ is invertible as a $(\C,\D)$-bimodule category if and only if
	\begin{subalign}[eqn:invertibilityConditions]
		 & \FPdim{\C} = \FPdim{\D} \hspace{10mm} \mathrm{and} \label{eqn:invertibilityDimensionCondition} \\
		 & \frac{1}{\rk\M}
		\sum_{\substack{\catOb[C]{a}\in\irr{\C}                                                           \\
		\catOb[M]{b},\catOb[M]{d} \in\irr{\M}                                                             \\
				\catOb[M]{\alpha},\catOb[M]{\beta},\catOb[M]{\mu},\catOb[M]{\nu}}}
		\frac{\dimC{a}}{\dimM{b}^2}\F[2]{a}{b}{c}{d}{\alpha,b,\mu}{\mu,d,\beta} \Fi[2]{a}{b}{c'}{d}{\alpha,b,\nu}{\nu,d,\beta}
		=\deltaD{c}{c^\prime},\label{eqn:annCatCharacterRelations}
	\end{subalign}
	for all $\catOb[D]{c},\,\catOb[D]{c^\prime}\in\irr{\D}$.
	\begin{proof}
		Given representations of $\ann{\C}{\M}$ labeled by $\catOb[D]{c},\,\catOb[D]{c^\prime}\in\irr{\D}$, \cref{eqn:annCatCharacterRelations} is $(\chi_{\catOb[D]{c}},\chi_{\catOb[D]{c^\prime}})$ evaluated using \cref{eqn:characterBasisVector}.

		`$\implies$': If $\M$ is invertible, then $\D\tensorEquiv\End{\C}{\M}$. In that case, simple objects of $\D$ label irreducible representations, so the associated irreducible characters are orthonormal by \cref{result:SchurOrthogonality}.

		`$\impliedby$': In the case $\M$ is not invertible, (simple) objects of $\D$ still label representations of $\ann{\C}{\M}$. The conditions \cref{eqn:invertibilityConditions} can fail in several ways:
		\begin{enumerate}
			\item Some irreducible representations are missing, for example is $\D$ is a subcategory of $\End{\C}{\M}$. In this case, \cref{eqn:annCatCharacterRelations} will still hold, but \cref{eqn:invertibilityDimensionCondition} will fail.

			\item Distinct objects of $\D$ label the same irreducible representation. In this case \cref{eqn:annCatCharacterRelations} will fail.

			\item Objects of $\D$ label reducible representations. In this case, \cref{eqn:invertibilityDimensionCondition} will fail since the questionable representation will have overlap $>1$ with itself.
		\end{enumerate}
		Multiple of these may occur for any given example.

		Finally, one could be concerned that objects of $\D$ label irreducible representations, but this is not a tensor equivalence. This is ruled out by the $\FName{3}$, which witnesses the tensor equivalence as discussed above.
	\end{proof}
\end{theorem}

All three failure modes of \cref{thm:invertibility} occur, for example
\begin{enumerate}
	\item The bimodule $\bimodule[e]{\vvec{\ZZ{2}}}{\vvec{}}{\vvec{}}$ is an example.
	\item Consider $\bimodule[e]{\vvec{\ZZ{2}}}{\vvec{}}{\vvec{\ZZ{2}}}$, with all bimodule associators equal to $+1$. In this case, both objects of $\D$ label the trivial representation of $\ZZ{2}$.
	\item The bimodule $\bimodule[e]{\vvec{\ZZ{2}}}{\vvec{}}{\rrep{\symmetricGroup{3}}}$ displays this failure. The 2-dimensional representation $\pi$ of $\symmetricGroup{3}$ labels the reducible representation $1\oplus\sigma$ of $\ZZ{2}$. This can be seen by restricting $\pi$ to $\ZZ{2}$. 
\end{enumerate}


\section{Schur orthogonality for invertible bimodules}\label{sec:orthogonality}

The results in \cref{sec:invertibility} are Schur's orthogonality relations for characters of WHAs. In this section, we discuss orthogonality relations for matrix elements of $\ann{\C}{\M}$. In addition to their intrinsic importance, these relations have application in physics. We discuss these applications in \cref{sec:MPOinjectivity}.
We note that these relations appeared in a different context in \onlinecite{Petkova2001}.

\begin{theorem}[Orthogonality of matrix elements]\label{thm:orthogMatrixElements}
	Let $\C$ be a unitary fusion category, and $\invertibleBimodule{\C}{\M}$ an indecomposable, unitary, finitely semisimple, invertible bimodule category.

	Let $\catOb[D]{c},\catOb[D]{c^\prime}$ be simple objects in $\End{\C}{\M}$, then
	\begin{align}
		\sum_{\substack{
		\catOb[C]{a} \\
				\catOb[M]{\alpha},\catOb[M]{\nu}}} \dimC{a} \F[2]{a}{b}{c}{d}{\alpha,e,\beta}{\mu,f,\nu} \Fi[2]{a}{b}{c'}{d}{\alpha,e,\beta '}{\mu ',f,\nu}=\deltaD{c}{c^\prime}\deltaM{\beta}{\beta^\prime}\deltaM{\mu}{\mu^\prime}  \frac{\dimM{e} \dimM{f}}{\dimD{c}}\label{eqn:matrixElementOrthog}
	\end{align}
\end{theorem}

To prove this, we will need several results. The first allows us to `change variables' in expressions involving the Haar integral\footnote{\cref{result:changeOfVariables} generalizes $\sum_g g^{-1}\otimes gx = \sum_g xg^{-1}\otimes g$, familiar from finite groups.}.
\begin{lemma}\label{result:changeOfVariables}

	Let $A$ be a WHA with Haar integral $\Lambda$, then
	\begin{align}
		S(\Lambda_{(1)})\otimes \Lambda_{(2)}x & =xS(\Lambda_{(1)})\otimes \Lambda_{(2)},\label{eqn:changeOfVariables}
	\end{align}
	for all $x\in A$.
	\begin{proof}
		Let $y=S^{-1}(x)$, then \cref{eqn:changeOfVariables} becomes
		\begin{subalign}
			S(\Lambda_{(1)})\otimes \Lambda_{(2)}S(y) & =S(y)S(\Lambda_{(1)})\otimes \Lambda_{(2)} \\
			                                          & = S(\Lambda_{(1)}y)\otimes \Lambda_{(2)}.
		\end{subalign}
		It is therefore sufficient to show that
		\begin{align}
			\Lambda_{(1)}\otimes \Lambda_{(2)}S(y) & = \Lambda_{(1)}y\otimes \Lambda_{(2)},
		\end{align}
		which is the right integral version of \onlinecite[Lemma 3.2]{Boehm1999}.
	\end{proof}
\end{lemma}

We also need to understand $A$-module homomorphisms.
\begin{lemma}\label{result:diagonalX}
	Let $V,W$ be simple $A$ modules, with action $\rho_V,\rho_W$. Let $M\in\Mat{\dim{W}\times\dim{V}}$. Define
	\begin{align}
		X_M := \rho_W\left(S(\Lambda_{(1)})\right) M \rho_V\left(\Lambda_{(2)}\right),\label{eqn:defXM}
	\end{align}
	then $X_M\in\Hom_A(V,W)$, and
	\begin{align}
		X_M = \delta_{V}^{W} \frac{\Tr{X_M}}{\dim{V}} \one_V.
	\end{align}

	\begin{proof}
		From \cref{result:changeOfVariables}, we have
		\begin{subalign}
			\rho_W\left(x\right) X_M & =\rho_W\left(x\right) \rho_W\left(S(\Lambda_{(1)})\right) M \rho_V\left(\Lambda_{(2)}\right)                                                                                                 \\
			                         & = \rho_W\left(x S(\Lambda_{(1)})\right) M \rho_V\left(\Lambda_{(2)}\right)  = \rho_W\left(S(\Lambda_{(1)})\right) M \rho_V\left(\Lambda_{(2)} x\right) \tag{\cref{result:changeOfVariables}} \\
			                         & = \rho_W\left(S(\Lambda_{(1)})\right) M \rho_V\left(\Lambda_{(2)}\right) \rho_V\left(x\right) = X_M \rho_V\left(x\right) ,
		\end{subalign}
		for any $x\in A$, so $X_M\in \Hom_A(V,W)$. By Schur's lemma, we know that
		\begin{align}
			X_M = \delta_{V}^{W} c_M \one_V, \label{eqn:schurUnknown}
		\end{align}
		for some $c_M\in\mathbb{C}$. Taking the trace of \cref{eqn:schurUnknown} completes the proof.
	\end{proof}
\end{lemma}
To understand the form of $c_M$, we need another result.
\begin{lemma}\label{result:integralToGroup}
	Let $V$ label a simple module, then
	\begin{align}
		\rho_V\left(\Lambda_{(2)}S\left(\Lambda_{(1)}\right)\right) & = \frac{\dim{V}}{\varepsilon(1)\dimQ{V}}\rho_V(g^{-1}),
	\end{align}
	where $g$ is the canonical grouplike element in the WHA~\cite{Boehm1999}.
	\begin{proof}
		Choose a set of matrix units $e^{r}_{ij}$ for the irreducible representations $r$ of $A$.
		\begin{subalign}
			S\left(\Lambda_{(1)}\right) \otimes \Lambda_{(2)} & = \sum_{\catOb[D]{r}\in \irr{\rrep{A}}} \frac{1}{\varepsilon(1)\dimD{r}}\sum_{i,j} e_{ij}^{r} g^{-1/2} \otimes g^{-1/2} e_{ji}^{r}   \tag{\cite[Lemma 4.6]{Boehm1999}} \\
			\implies \Lambda_{(2)}S\left(\Lambda_{(1)}\right) & = \sum_{\catOb[D]{r}\in \irr{\rrep{A}}} \frac{1}{\varepsilon(1)\dimD{r}}\sum_{i,j} g^{-1/2} e_{ji}^{r}e_{ij}^{r} g^{-1/2}                                              \\
			                                                  & = \sum_{\catOb[D]{r}\in \irr{\rrep{A}}} \frac{1}{\varepsilon(1)\dimD{r}}\sum_{i,j} g^{-1/2} e_{jj}^{r} g^{-1/2}                                                        \\
			                                                  & = \sum_{\catOb[D]{r}\in \irr{\rrep{A}}} \frac{\dim{V_{\catOb[D]{r}}}}{\varepsilon(1)\dimD{r}}g^{-1/2} \one_{\catOb[D]{r}} g^{-1/2}.
		\end{subalign}

	\end{proof}
\end{lemma}

Using \cref{result:integralToGroup}, we have
\begin{subalign}
	c_M & = \frac{1}{\varepsilon(1)\dimQ{V}} \Tr\left(\rho_V(g^{-1}) M\right),                      \\
	X_M & = \delta_{V}^{W} \frac{\Tr\left(\rho_V(g^{-1}) M\right)}{\varepsilon(1)\dimQ{V}}  \one_V.
\end{subalign}

Finally, we can prove \cref{thm:orthogMatrixElements}.
\begin{proof}[Proof of \cref{thm:orthogMatrixElements}]
	In \cref{eqn:action}, we introduced a representation of $\ann{\C}{\M}$.
	In this representation, the grouplike element (see \cref{eqn:grouplike}) has matrix entries
	\begin{align}
		\rho_V(g^{-1})_{(\printColorList{k_0,\kappa,k_1}{MColor,MColor,MColor}),(\printColorList{j_0,\varphi,j_1}{MColor,MColor,MColor})} & = \deltaM{j_{0}}{k_{0}} \deltaM{j_{1}}{k_{1}} \deltaM{\kappa}{\varphi} \frac{ \dimM{j_{1}}}{\dimM{j_{0}}}.
	\end{align}

	Specializing to $M= e_{(j_0,\varphi,j_1),(k_0,\kappa,k_1)} $, we have
	\begin{subalign}
		X_M & = \deltaD{c}{c^{\prime}} \frac{1}{\varepsilon(1)\dimD{c}}  \deltaM{j_{0}}{k_{0}} \deltaM{j_{1}}{k_{1}} \deltaM{\kappa}{\varphi} \frac{ \dimM{j_{1}}}{\dimM{j_{0}}} \one_{V_{\catOb[D]{c}}} \\
		    & = \deltaD{c}{c^{\prime}} \frac{1}{\rk\M\dimD{c}} \deltaM{j_{0}}{k_{0}} \deltaM{j_{1}}{k_{1}} \deltaM{\kappa}{\varphi} \frac{ \dimM{j_{1}}}{\dimM{j_{0}}} \one_{V_{\catOb[D]{c}}}.
	\end{subalign}

	The result holds by computing matrix elements of $X_M$ on both sides of \cref{eqn:defXM}.
\end{proof}

\section{Applications}\label{sec:applications}

We now discuss two applications of the results in \cref{sec:invertibility,sec:orthogonality}, the first concerns fusion category data and the second is a physical application.


\subsection{Computing data}\label{sec:computing}

In this section, we discuss how the representation theoretic definition of the dual category allows computation of all skeletal data for an invertible bimodule category. This extends the work initiated in \onlinecites{Bridgeman2020a,Barter2022}, where only the computation of $\FName{4}$ was considered.

As input, we require $(\FName{0},\FName{1})$, the skeletal data for a unitary module category $\leftmodule{\C}{\M}$. After running the algorithm we describe, one obtains $(\FName{2},\FName{3},\FName{4})$ for the invertible bimodule category $\invertibleBimodule{\C}{\M}$, in a unitary gauge. This data is crucial, for example, in the tensor networks discussed in \onlinecite{Lootens2021}, and in \cref{sec:MPOinjectivity}, as it is the input data to those constructions.

The algorithm proceeds as follows:
\begin{enumerate}
	\item Compute all irreducible representations of $\ann{\C}{\M}$, with an orthonormal basis.
	\item Compute $\FName{2}$ using \cref{eqn:action}.
	\item Compute $\FName{3}$ as trivalent intertwining maps.
	\item Compute $\FName{4}$ as discussed in \onlinecite{Barter2022}.
\end{enumerate}

The irreducible representations can be computed as discussed in \onlinecite{Barter2022}. Here, unlike the previous work, we insist the basis for each representation is chosen to be orthonormal. Additionally, we choose a $*$-representation, so $\langle v, T w\rangle = \langle T^* v,w\rangle$, for each annular diagram $T$ and pair of vectors $v,\,w$ in the representation.
With these bases fixed, $\FName{2}$ is given by the structure constants of the representation via \cref{eqn:action}. 

Given explicit representations $a,\,b,\,c$, intertwining maps $V_{ab}^{c;\alpha}$, which embed $c$ into $a\otimes b$, can be chosen (see \onlinecite{Barter2022}). To ensure unitarity, they should be chosen to embed the representations isometrically. The matrix elements of these maps determine $\FName{3}$:
\begin{align}
	\Fi[3]{a}{b}{c}{d}{\alpha,e,\beta}{\mu,f,\nu} &= \left[ V_{bc}^{f;\mu}  \right]_{\hat{I},\hat{J}},
\end{align}
where
\begin{align}
	\hat{I} = \frac{I}{\norm{I}}&&\hat{J} = \frac{J}{\norm{J}}&&
	I = \begin{tikzarray}[scale=.3]{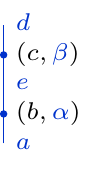}
		\begin{scope}[shift={(0,-1)}]
			\irrepVector{a}{e}{b}{\alpha}
		\end{scope}
		\begin{scope}[shift={(0,1)}]
			\irrepVector{}{d}{c}{\beta}
		\end{scope}
	\end{tikzarray}
	&&J = \begin{tikzarray}[scale=.3]{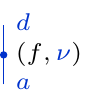}
			\irrepVector{a}{d}{f}{\nu}
	\end{tikzarray}.
\end{align}
Finally, $\FName{4}$ can be computed by relating these intertwiners as discussed in \onlinecite{Barter2022}. Isometricity of the intertwining maps $V$ fixes unitarity of $\FName{4}$~\cite{Barter2022}. Expressing the $V$'s with respect to an orthonormal basis ensures unitarity of $\FName{3}$, while unitarity of $\FName{2}$ is fixed by choosing the representations to be compatible with the $*$-structure on $\ann{\C}{\M}$:

\subsubsection*{Unitarity of \texorpdfstring{${}^{\triangleright\hspace{-.1em}\triangleleft\!}F$}{bimodule associator}}

Let $\catOb[D]{c}$ label some irreducible representation of $\Tub{\C}{\M}$, and define the inner product on $V_{\catOb[D]{c}}$ so that the `dot basis' (left side in \cref{eqn:relateBases}) is orthonormal, then for the action of annular diagrams to respect the $*$-structure
\begin{subalign}
	\sum_{
		\catOb[M]{e},\catOb[M]{\alpha}}\sqrt{\dimM{e}} &
	\begin{tikzarray}[scale=.55]{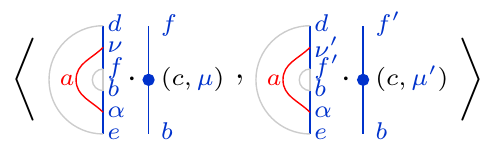}
		\begin{scope}[local bounding box=sc1]
			\node at (0,0) {$\biggl\langle$};
		\end{scope}
		\begin{scope}[local bounding box=sc2,anchor=west,shift={($(sc1.east)+(1,0)$)}]
			\tubd{b,f}{e,d}{\alpha,a,\nu}
		\end{scope}
		\begin{scope}[local bounding box=sc3,anchor=west,shift={($(sc2.east)+(-.25,0)$)}]
			\node at (0,0) {$\cdot$};
		\end{scope}
		\begin{scope}[local bounding box=sc4,anchor=west,shift={($(sc3.east)+(0,0)$)}]
			\irrepVector{b}{f}{c}{\mu};
		\end{scope}
		\begin{scope}[local bounding box=sc5,anchor=west,shift={($(sc4.east)+(-.25,0)$)}]
			\node at (0,0) {$,$};
		\end{scope}
		\begin{scope}[local bounding box=sc6,anchor=west,shift={($(sc5.east)+(1,0)$)}]
			\tubd{b,f'}{e,d}{\alpha,a,\nu '};
		\end{scope}
		\begin{scope}[local bounding box=sc7,anchor=west,shift={($(sc6.east)+(-.25,0)$)}]
			\node at (0,0) {$\cdot$};
		\end{scope}
		\begin{scope}[local bounding box=sc8,anchor=west,shift={($(sc7.east)+(0,0)$)}]
			\irrepVector{b}{f'}{c}{\mu '};
		\end{scope}
		\begin{scope}[local bounding box=sc9,anchor=west,shift={($(sc8.east)+(-.25,0)$)}]
			\node at (0,0) {$\biggr\rangle$};
		\end{scope}
	\end{tikzarray}
	= \sum_{
		\catOb[M]{e},\catOb[M]{\alpha}}\sqrt{\dimM{e}}
	\begin{tikzarray}[scale=.55]{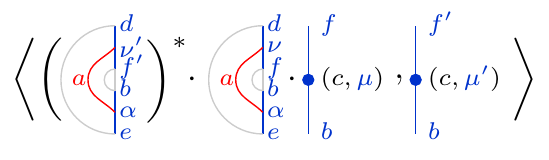}
		\begin{scope}[local bounding box=sc1]
			\node at (0,0) {$\biggl\langle\biggl($};
		\end{scope}
		\begin{scope}[local bounding box=sc2,anchor=west,shift={($(sc1.east)+(.75,0)$)}]
			\tubd{b,f'}{e,d}{\alpha,a,\nu '}
		\end{scope}
		\begin{scope}[local bounding box=sc3,anchor=west,shift={($(sc2.east)+(-.25,0)$)}]
			\node at (0,0) {$\biggr)^{*}$};
			\node at (.75,0) {$\cdot$};
		\end{scope}
		\begin{scope}[local bounding box=sc4,anchor=west,shift={($(sc3.east)+(1,0)$)}]
			\tubd{b,f}{e,d}{\alpha,a,\nu}
		\end{scope}
		\begin{scope}[local bounding box=sc5,anchor=west,shift={($(sc4.east)+(-.25,0)$)}]
			\node at (0,0) {$\cdot$};
		\end{scope}
		\begin{scope}[local bounding box=sc6,anchor=west,shift={($(sc5.east)+(0,0)$)}]
			\irrepVector{b}{f}{c}{\mu};
		\end{scope}
		\begin{scope}[local bounding box=sc7,anchor=west,shift={($(sc6.east)+(-.25,0)$)}]
			\node at (0,0) {$,$};
		\end{scope}
		\begin{scope}[local bounding box=sc8,anchor=west,shift={($(sc7.east)+(0,0)$)}]
			\irrepVector{b}{f'}{c}{\mu '};
		\end{scope}
		\begin{scope}[local bounding box=sc9,anchor=west,shift={($(sc8.east)+(-.25,0)$)}]
			\node at (0,0) {$\biggr\rangle$};
		\end{scope}
	\end{tikzarray}
	\\
	& \implies
	\sum_{\catOb[M]{e},\catOb[M]{\alpha},\catOb[M]{\beta}} \F[2]{a}{b}{c}{d}{\alpha,e,\beta}{\mu,f,\nu} \Fi[2]{a}{b}{c}{d}{\alpha,e,\beta}{\mu',f',\nu'}
	\sqrt{\frac{\dimM{b} \dimM{f}}{\dimM{d} \dimM{e}}}
	\frac{\X{b}{c}{f}{\mu} \X{b}{c}{f}{\mu}^*}{\X{e}{c}{d}{\beta} \X{e}{c}{d}{\beta}^*}
	= \deltaM{f}{f'} \deltaM{\mu}{\mu '} \deltaM{\nu}{\nu '}.
\end{subalign}
Choosing
\begin{align}
	\X{a}{b}{c}{\alpha} & = \frac{\tau_{\catOb[D]{b}}}{(\dimM{a}\dimM{c})^{1/4}},
\end{align}
where $\tau_{\catOb[D]{b}}$ is arbitrary, means that $\FName{2}$ is unitary, and the representation is a $*$-representation.

In \onlinecite{VecGCode}, we provide a Mathematica package for computing the data in the case $\C=\vvec{G}$, for finite groups $G$.

\subsection{MPO-injectivity}\label{sec:MPOinjectivity}

We now turn to an application of the generalized Schur orthogonality condition to the theory of topological order and tensor networks. In (2+1)-dimensions, it has been understood that tensor network representations of topologically ordered ground states, represented as projected entangled pair states (PEPS), satisfy a condition called MPO-injectivity \cite{Sahinoglu2021}. Generalizing the notion of $G$-injectivity~\cite{Schuch2011},
it states that when interpreted as a map from the virtual to the physical space, the PEPS tensors are injective on a subspace of the virtual, and that the projector onto this subspace is given in terms of a projector matrix product operator (MPO). A key feature of MPO-injective PEPS is that the ground state degeneracy does not grow with the system size. This is a necessary condition for topologically ordered fixed point models where the degeneracy should only depend on the topology of the underlying manifold. Additionally, these topologically ordered PEPS exhibit virtual MPO symmetries; string-like operators that can be moved freely through the lattice. The ground states spanned by the PEPS tensors can be characterized using these virtual MPO symmetries via their Ocneanu tube algebra~\cite{Ocneanu1993}, and MPO-injectivity guarantees that the number of distinct ground states is upper bounded by the number of irreps of this algebra.

In \onlinecite{Bultinck2017}, it was shown that such projector MPOs are characterized by fusion categories, and in \onlinecite{Lootens2021}, a general description of PEPS representations of string-net models and their MPO symmetries was provided in terms of the skeletal data of bimodule categories. For a string-net model based on an input fusion category $\D$, PEPS realizations of the ground state are determined by a choice of (right) module category $\rightmodule{\M}{\D}$ while their symmetries are described by a fusion category $\C$ such that $\bimodule{\C}{\M}{\D}$ is a $(\C,\D)$-bimodule category. Explicitly, the nonzero components of the PEPS tensors evaluate to
\begin{align}\label{eq:PEPSdef}
	\begin{tikzarray}[scale=.15, every node/.style={scale=0.8}]{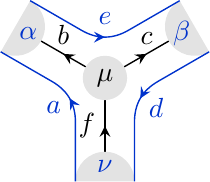}
		\begin{scope}
			\draw [draw=DColor,mid arrow] (0,-5) -- (0,-1.5) node[midway,left] {$\catOb[D]{f}$};
			\fill[gray!22] (-2,-7) -- (2,-7) arc (0:180:2) -- cycle;
			\coordinate (l11) at (-2,-7);
			\coordinate (l21) at (-2,-3);
			\coordinate (r11) at (2,-7);
			\coordinate (r21) at (2,-3);
			\node at (0,-6) {$\catOb[M]{\nu}$};
			\node at (0,4) {$\catOb[M]{e}$};
		\end{scope}
		\begin{scope}[rotate=120]
			\draw [draw=DColor,mid arrow] (0,-1.5) -- (0,-5) node[midway,above] {$\catOb[D]{c}$};
			\fill[gray!22] (-2,-7) -- (2,-7) arc (0:180:2) -- cycle;
			\coordinate (l12) at (-2,-7);
			\coordinate (l22) at (-2,-3);
			\coordinate (r12) at (2,-7);
			\coordinate (r22) at (2,-3);
			\node at (0,-6) {$\catOb[M]{\beta}$};
			\node at (0,4) {$\catOb[M]{a}$};
		\end{scope}
		\begin{scope}[rotate=240]
			\draw [draw=DColor,mid arrow] (0,-1.5) -- (0,-5) node[midway,above] {$\catOb[D]{b}$};
			\fill[gray!22] (-2,-7) -- (2,-7) arc (0:180:2) -- cycle;
			\coordinate (l13) at (-2,-7);
			\coordinate (l23) at (-2,-3);
			\coordinate (r13) at (2,-7);
			\coordinate (r23) at (2,-3);
			\node at (0,-6) {$\catOb[M]{\alpha}$};
			\node at (0,4) {$\catOb[M]{d}$};
		\end{scope}
		\draw [MColor,mid arrow] (l11) -- (l21) to [out=90,in=-30] (r23) -- (r13);
		\draw [MColor,mid arrow] (l12) -- (l22) to [out=210,in=90] (r21) -- (r11);
		\draw [MColor,mid arrow] (l13) -- (l23) to [out=330,in=210] (r22) -- (r12);
		\fill[gray!22] (0,0) circle (1.5);
		\node at (0,0) {$\catOb[D]{\mu}$};
	\end{tikzarray} = \F[3]{a}{b}{c}{d}{\alpha,e,\beta}{\mu,f,\nu} =: \hspace{-2em}
	\begin{tikzarray}[scale=.15, every node/.style={scale=0.8}]{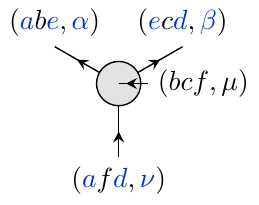}
		\begin{scope}
			\draw [draw=DColor,mid arrow] (0,-5) -- (0,-1.5) node[at start,below] {$(\catOb[M]{a}\catOb[D]{f}\catOb[M]{d},\catOb[M]{\nu})$};
		\end{scope}
		\begin{scope}[rotate=120]
			\draw [draw=DColor,mid arrow] (0,-1.5) -- (0,-5) node[above] {$(\catOb[M]{e}\catOb[D]{c}\catOb[M]{d},\catOb[M]{\beta})$};
		\end{scope}
		\begin{scope}[rotate=240]
			\draw [draw=DColor,mid arrow] (0,-1.5) -- (0,-5) node[above] {$(\catOb[M]{a}\catOb[D]{b}\catOb[M]{e},\catOb[M]{\alpha})$};
		\end{scope}
		\draw[fill=gray!22] (0,0) circle (1.5);
		\draw [draw=DColor,Q arrow reversed] (0,0) -- (2,0) node[right] {$(\catOb[D]{b}\catOb[D]{c}\catOb[D]{f},\catOb[D]{\mu})$};
	\end{tikzarray}
	, \quad
	\begin{tikzarray}[yscale=-.15,xscale=.15, every node/.style={scale=0.8}]{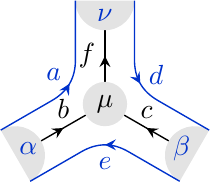}
		\begin{scope}
			\draw [draw=DColor,mid arrow reversed] (0,-5) -- (0,-1.5) node[midway,left] {$\catOb[D]{f}$};
			\fill[gray!22] (-2,-7) -- (2,-7) arc (0:180:2) -- cycle;
			\coordinate (l11) at (-2,-7);
			\coordinate (l21) at (-2,-3);
			\coordinate (r11) at (2,-7);
			\coordinate (r21) at (2,-3);
			\node at (0,-6) {$\catOb[M]{\nu}$};
			\node at (0,4) {$\catOb[M]{e}$};
		\end{scope}
		\begin{scope}[rotate=120]
			\draw [draw=DColor,mid arrow reversed] (0,-1.5) -- (0,-5) node[midway,above] {$\catOb[D]{c}$};
			\fill[gray!22] (-2,-7) -- (2,-7) arc (0:180:2) -- cycle;
			\coordinate (l12) at (-2,-7);
			\coordinate (l22) at (-2,-3);
			\coordinate (r12) at (2,-7);
			\coordinate (r22) at (2,-3);
			\node at (0,-6) {$\catOb[M]{\beta}$};
			\node at (0,4) {$\catOb[M]{a}$};
		\end{scope}
		\begin{scope}[rotate=240]
			\draw [draw=DColor,mid arrow reversed] (0,-1.5) -- (0,-5) node[midway,above] {$\catOb[D]{b}$};
			\fill[gray!22] (-2,-7) -- (2,-7) arc (0:180:2) -- cycle;
			\coordinate (l13) at (-2,-7);
			\coordinate (l23) at (-2,-3);
			\coordinate (r13) at (2,-7);
			\coordinate (r23) at (2,-3);
			\node at (0,-6) {$\catOb[M]{\alpha}$};
			\node at (0,4) {$\catOb[M]{d}$};
		\end{scope}
		\draw [MColor,mid arrow reversed] (l11) -- (l21) to [out=90,in=-30] (r23) -- (r13);
		\draw [MColor,mid arrow reversed] (l12) -- (l22) to [out=210,in=90] (r21) -- (r11);
		\draw [MColor,mid arrow reversed] (l13) -- (l23) to [out=330,in=210] (r22) -- (r12);
		\fill[gray!22] (0,0) circle (1.5);
		\node at (0,0) {$\catOb[D]{\mu}$};
	\end{tikzarray} = \Fi[3]{a}{b}{c}{d}{\alpha,e,\beta}{\mu,f,\nu} := \hspace{-2em}
	\begin{tikzarray}[yscale=-.15,xscale=.15, every node/.style={scale=0.8}]{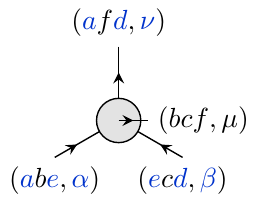}
		\begin{scope}
			\draw [draw=DColor,mid arrow reversed] (0,-5) -- (0,-1.5) node[at start,above] {$(\catOb[M]{a}\catOb[D]{f}\catOb[M]{d},\catOb[M]{\nu})$};
		\end{scope}
		\begin{scope}[rotate=120]
			\draw [draw=DColor,mid arrow reversed] (0,-1.5) -- (0,-5) node[below] {$(\catOb[M]{e}\catOb[D]{c}\catOb[M]{d},\catOb[M]{\beta})$};
		\end{scope}
		\begin{scope}[rotate=240]
			\draw [draw=DColor,mid arrow reversed] (0,-1.5) -- (0,-5) node[below] {$(\catOb[M]{a}\catOb[D]{b}\catOb[M]{e},\catOb[M]{\alpha})$};
		\end{scope}
		\draw[fill=gray!22] (0,0) circle (1.5);
		\draw [draw=DColor,mid arrow] (0,0) -- (2,0) node[right] {$(\catOb[D]{b}\catOb[D]{c}\catOb[D]{f},\catOb[D]{\mu})$};
	\end{tikzarray},
\end{align}
while the nonzero components of the MPO tensors are given by
\begin{align}
	\begin{tikzarray}[scale=.15, every node/.style={scale=0.8}]{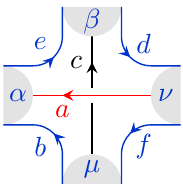}
		\def\ds{1}
		\begin{scope}
			\draw [DColor] (0,-.5) -- (0,-5+\ds);
			\fill[gray!22] (-2,-7+\ds) -- (2,-7+\ds) arc (0:180:2) -- cycle;
			\coordinate (l11) at (-2,-7+\ds);
			\coordinate (l21) at (-2,-4);
			\coordinate (r11) at (2,-7+\ds);
			\coordinate (r21) at (2,-4);
			\node at (-3.5,-3.5) {$\catOb[M]{b}$};
			\node at (0,-6+\ds) {$\catOb[M]{\mu}$};
		\end{scope}
		\begin{scope}[rotate=90]
			\draw [CColor] (0,0) -- (0,-5+\ds);
			\fill[gray!22] (-2,-7+\ds) -- (2,-7+\ds) arc (0:180:2) -- cycle;
			\coordinate (l12) at (-2,-7+\ds);
			\coordinate (l22) at (-2,-4);
			\coordinate (r12) at (2,-7+\ds);
			\coordinate (r22) at (2,-4);
			\node at (-3.5,-3.5) {$\catOb[M]{f}$};
			\node at (0,-6+\ds) {$\catOb[M]{\nu}$};
		\end{scope}
		\begin{scope}[rotate=180]
			\draw [draw=DColor,mid arrow] (0,-.5) -- (0,-5+\ds) node[midway,left] {$\catOb[D]{c}$};
			\fill[gray!22] (-2,-7+\ds) -- (2,-7+\ds) arc (0:180:2) -- cycle;
			\coordinate (l13) at (-2,-7+\ds);
			\coordinate (l23) at (-2,-4);
			\coordinate (r13) at (2,-7+\ds);
			\coordinate (r23) at (2,-4);
			\node at (-3.5,-3.5) {$\catOb[M]{d}$};
			\node at (0,-6+\ds) {$\catOb[M]{\beta}$};
		\end{scope}
		\begin{scope}[rotate=270]
			\draw [draw=CColor,mid arrow] (0,0) -- (0,-5+\ds) node[midway,below] {$\catOb[C]{a}$};;
			\fill[gray!22] (-2,-7+\ds) -- (2,-7+\ds) arc (0:180:2) -- cycle;
			\coordinate (l14) at (-2,-7+\ds);
			\coordinate (l24) at (-2,-4);
			\coordinate (r14) at (2,-7+\ds);
			\coordinate (r24) at (2,-4);
			\node at (-3.5,-3.5) {$\catOb[M]{e}$};
			\node at (0,-6+\ds) {$\catOb[M]{\alpha}$};
		\end{scope}
		\draw [MColor,mid arrow] (l11) -- (l21) to [out=90,in=0] (r24) -- (r14);
		\draw [MColor,mid arrow] (l12) -- (l22) to [out=180,in=90] (r21) -- (r11);
		\draw [MColor,mid arrow] (l13) -- (l23) to [out=270,in=180] (r22) -- (r12);
		\draw [MColor,mid arrow] (l14) -- (l24) to [out=0,in=270] (r23) -- (r13);
	\end{tikzarray} = \F[2]{a}{b}{c}{d}{\alpha,e,\beta}{\mu,f,\nu} =:
	\begin{tikzarray}[scale=.15, every node/.style={scale=0.8}]{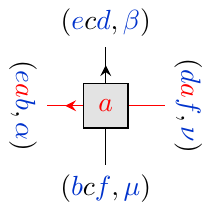}
		\def\ds{1}
		\begin{scope}
			\draw [draw=DColor] (0,-1.5) -- (0,-5+\ds) node[below] {$(\catOb[M]{b}\catOb[D]{c}\catOb[M]{f},\catOb[M]{\mu})$};
		\end{scope}
		\begin{scope}[rotate=90]
			\draw [draw=CColor] (0,-1.5) -- (0,-5+\ds) node[right] {\rotatebox{-90}{$(\catOb[M]{d}\catOb[C]{a}\catOb[M]{f},\catOb[M]{\nu})$}};
		\end{scope}
		\begin{scope}[rotate=180]
			\draw [draw=DColor,mid arrow] (0,-1.5) -- (0,-5+\ds) node[above] {$(\catOb[M]{e}\catOb[D]{c}\catOb[M]{d},\catOb[M]{\beta})$};
		\end{scope}
		\begin{scope}[rotate=270]
			\draw [draw=CColor,mid arrow] (0,-1.5) -- (0,-5+\ds) node[left] {\rotatebox{-90}{$(\catOb[M]{e}\catOb[C]{a}\catOb[M]{b},\catOb[M]{\alpha})$}};
		\end{scope}
		\draw[fill=gray!22] (-1.5,-1.5) rectangle (1.5,1.5);
		\node[] at (0,0) {$\catOb[C]{a}$}
	\end{tikzarray}
	, \quad
	\begin{tikzarray}[xscale=-.15,yscale=.15, every node/.style={scale=0.8}]{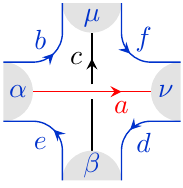}
		\def\ds{1}
		\begin{scope}
			\draw [DColor] (0,-.5) -- (0,-5+\ds);
			\fill[gray!22] (-2,-7+\ds) -- (2,-7+\ds) arc (0:180:2) -- cycle;
			\coordinate (l11) at (-2,-7+\ds);
			\coordinate (l21) at (-2,-4);
			\coordinate (r11) at (2,-7+\ds);
			\coordinate (r21) at (2,-4);
			\node at (-3.5,-3.5) {$\catOb[M]{d}$};
			\node at (0,-6+\ds) {$\catOb[M]{\beta}$};
		\end{scope}
		\begin{scope}[rotate=90]
			\draw [CColor] (0,0) -- (0,-5+\ds);
			\fill[gray!22] (-2,-7+\ds) -- (2,-7+\ds) arc (0:180:2) -- cycle;
			\coordinate (l12) at (-2,-7+\ds);
			\coordinate (l22) at (-2,-4);
			\coordinate (r12) at (2,-7+\ds);
			\coordinate (r22) at (2,-4);
			\node at (-3.5,-3.5) {$\catOb[M]{e}$};
			\node at (0,-6+\ds) {$\catOb[M]{\alpha}$};
		\end{scope}
		\begin{scope}[rotate=180]
			\draw [DColor,mid arrow] (0,-.5) -- (0,-5+\ds) node[midway,left] {$\catOb[D]{c}$};
			\fill[gray!22] (-2,-7+\ds) -- (2,-7+\ds) arc (0:180:2) -- cycle;
			\coordinate (l13) at (-2,-7+\ds);
			\coordinate (l23) at (-2,-4);
			\coordinate (r13) at (2,-7+\ds);
			\coordinate (r23) at (2,-4);
			\node at (-3.5,-3.5) {$\catOb[M]{b}$};
			\node at (0,-6+\ds) {$\catOb[M]{\mu}$};
		\end{scope}
		\begin{scope}[rotate=270]
			\draw [CColor,mid arrow] (0,0) -- (0,-5+\ds) node[midway,below] {$\catOb[C]{a}$};;
			\fill[gray!22] (-2,-7+\ds) -- (2,-7+\ds) arc (0:180:2) -- cycle;
			\coordinate (l14) at (-2,-7+\ds);
			\coordinate (l24) at (-2,-4);
			\coordinate (r14) at (2,-7+\ds);
			\coordinate (r24) at (2,-4);
			\node at (-3.5,-3.5) {$\catOb[M]{f}$};
			\node at (0,-6+\ds) {$\catOb[M]{\nu}$};
		\end{scope}
		\draw [MColor,mid arrow reversed] (l11) -- (l21) to [out=90,in=0] (r24) -- (r14);
		\draw [MColor,mid arrow reversed] (l12) -- (l22) to [out=180,in=90] (r21) -- (r11);
		\draw [MColor,mid arrow reversed] (l13) -- (l23) to [out=270,in=180] (r22) -- (r12);
		\draw [MColor,mid arrow reversed] (l14) -- (l24) to [out=0,in=270] (r23) -- (r13);
	\end{tikzarray} = \Fi[2]{a}{b}{c}{d}{\alpha,e,\beta}{\mu,f,\nu} =:
	\begin{tikzarray}[xscale=-.15,yscale=.15, every node/.style={scale=0.8}]{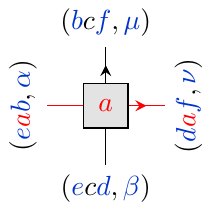}
		\def\ds{1}
		\begin{scope}
			\draw [DColor] (0,-1.5) -- (0,-5+\ds) node[below] {$(\catOb[M]{e}\catOb[D]{c}\catOb[M]{d},\catOb[M]{\beta})$};
		\end{scope}
		\begin{scope}[rotate=90]
			\draw [draw=CColor] (0,-1.5) -- (0,-5+\ds) node[left] {\rotatebox{90}{$(\catOb[M]{e}\catOb[C]{a}\catOb[M]{b},\catOb[M]{\alpha})$}};
		\end{scope}
		\begin{scope}[rotate=180]
			\draw [DColor,mid arrow] (0,-1.5) -- (0,-5+\ds) node[above] {$(\catOb[M]{b}\catOb[D]{c}\catOb[M]{f},\catOb[M]{\mu})$};
		\end{scope}
		\begin{scope}[rotate=270]
			\draw [draw=CColor,mid arrow] (0,-1.5) -- (0,-5+\ds) node[right] {\rotatebox{90}{$(\catOb[M]{d}\catOb[C]{a}\catOb[M]{f},\catOb[M]{\nu})$}};
		\end{scope}
		\draw[fill=gray!22] (-1.5,-1.5) rectangle (1.5,1.5);
		\node at (0,0) {$\catOb[C]{a}$};
	\end{tikzarray}
	.
\end{align}
It is well known that the distinct ground-states of a string-net model $\D$ are in correspondence with the monoidal center $\drinfeld{\D}$. On the other hand, as mentioned above, the distinct ground states are also characterized as irreps of the MPO symmetry annular algebra, which are equivalent to $\drinfeld{\C}$. In order for these two characterizations to coincide, we need $\C$ and $\D$ to be Morita equivalent. We will now show that the Morita equivalence follows from imposing MPO-injectivity. Using the diagrammatic notation for tensor networks, the MPO-injectivity condition can graphically be depicted as\footnote{The factors $\dimM{e},\dimM{f}$ are added in accordance to the convention where we insert quantum dimensions for closed loops \cite{Williamson2017}, and the factor $\dimD{c'}$ is a normalization factor of the PEPS tensors that ensures the LHS is a projector.}
\begin{align}
	\frac{\dimM{e}\dimM{f}}{\dimD{c'}\FPdim{\D}}
	\begin{tikzarray}[scale=.15, every node/.style={scale=0.8}]{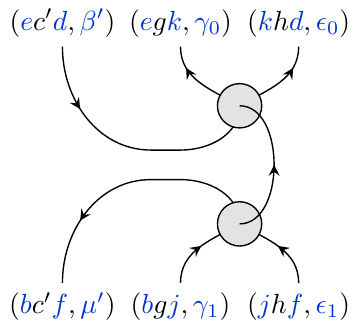}
		\draw [draw=DColor,TQ arrow] (0,4) to [out=30,in=270] (4,8);
		\draw [draw=DColor,TQ arrow] (0,4) to [out=150,in=270] (-4,8);
		\draw [draw=DColor,TQ arrow reversed] (0,4) to [out=270,in=0] (-4,1) -- (-6,1) to [out=180,in=270] (-12,8);
		\draw [draw=DColor,mid arrow reversed] (0,-4) to [out=-30,in=-270] (4,-8);
		\draw [draw=DColor,mid arrow reversed] (0,-4) to [out=-150,in=-270] (-4,-8);
		\draw [draw=DColor,TQ arrow] (0,-4) to [out=90,in=0] (-4,-1) -- (-6,-1) to [out=180,in=90] (-12,-8);
		\draw[fill=gray!22] (0,4) circle (1.5);
		\draw[fill=gray!22] (0,-4) circle (1.5);
		\draw [draw=DColor,mid arrow reversed] (0,4) to [out=0,in=0] (0,-4);
		\node[above] at (-12,8) {$(\catOb[M]{e}\catOb[D]{c'}\catOb[M]{d},\catOb[M]{\beta'})$};
		\node[above] at (-4,8) {$(\catOb[M]{e}\catOb[D]{g}\catOb[M]{k},\catOb[M]{\gamma_{0}})$};
		\node[above] at (4,8) {$(\catOb[M]{k}\catOb[D]{h}\catOb[M]{d},\catOb[M]{\epsilon_{0}})$};
		\node[below] at (-12,-8) {$(\catOb[M]{b}\catOb[D]{c'}\catOb[M]{f},\catOb[M]{\mu'})$};
		\node[below] at (-4,-8) {$(\catOb[M]{b}\catOb[D]{g}\catOb[M]{j},\catOb[M]{\gamma_{1}})$};
		\node[below] at (4,-8) {$(\catOb[M]{j}\catOb[D]{h}\catOb[M]{f},\catOb[M]{\epsilon_{1}})$};
		%
	\end{tikzarray} = \sum_{\catOb[C]{a}} \frac{\dimC{a}}{\FPdim{\C}}
	\begin{tikzarray}[scale=.15, every node/.style={scale=0.8}]{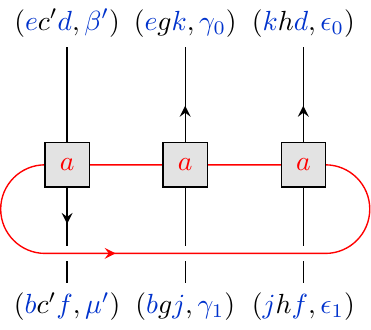}
		\draw [DColor,TQ arrow] (0,-8) -- (0,8);
		\draw [DColor,Q arrow reversed] (-8,-8) -- (-8,8);
		\draw [DColor,TQ arrow] (8,-8) -- (8,8);
		\fill[white] (-0.5,-6.5) rectangle (0.5,-5.5);
		\fill[white] (-8.5,-6.5) rectangle (-7.5,-5.5);
		\fill[white] (7.5,-6.5) rectangle (8.5,-5.5);
		\draw[CColor,TQ arrow reversed] (-9.5,0) -- (9.5,0) arc (90:-90:3) -- (-9.5,-6) arc (270:90:3);
		\draw[fill=gray!22] (-1.5,-1.5) rectangle (1.5,1.5);
		\draw[fill=gray!22] (-9.5,-1.5) rectangle (-6.5,1.5);
		\draw[fill=gray!22] (6.5,-1.5) rectangle (9.5,1.5);
		\node at (0,0) {$\catOb[C]{a}$};
		\node at (8,0) {$\catOb[C]{a}$};
		\node at (-8,0) {$\catOb[C]{a}$};
		\node[above] at (-8,8) {$(\catOb[M]{e}\catOb[D]{c'}\catOb[M]{d},\catOb[M]{\beta'})$};
		\node[above] at (-0,8) {$(\catOb[M]{e}\catOb[D]{g}\catOb[M]{k},\catOb[M]{\gamma_{0}})$};
		\node[above] at (8,8) {$(\catOb[M]{k}\catOb[D]{h}\catOb[M]{d},\catOb[M]{\epsilon_{0}})$};
		\node[below] at (-8,-8) {$(\catOb[M]{b}\catOb[D]{c'}\catOb[M]{f},\catOb[M]{\mu'})$};
		\node[below] at (0,-8) {$(\catOb[M]{b}\catOb[D]{g}\catOb[M]{j},\catOb[M]{\gamma_{1}})$};
		\node[below] at (8,-8) {$(\catOb[M]{j}\catOb[D]{h}\catOb[M]{f},\catOb[M]{\epsilon_{1}})$};
		%
	\end{tikzarray}
\end{align}
which when written out explicitly becomes
\begin{align}
	\frac{\dimM{e} \dimM{f}}{\dimD{c'} \FPdim{\D}}\sum_{\catOb[D]{\zeta}}  \Fi[3]{b}{g}{h}{f}{\gamma_{1},j,\epsilon_{1}}{\zeta,c',\mu '}\F[3]{e}{g}{h}{d}{\gamma_{0},k,\epsilon_{0}}{\zeta,c',\beta '} & =
	\sum_{\substack{\catOb[C]{a}                                                                                                                                                                           \\\catOb[M]{\alpha},\catOb[M]{\nu},\catOb[M]{\eta}}} \frac{\dimC{a}}{\FPdim{\C}}\Fi[2]{a}{b}{c'}{d}{\alpha,e,\beta '}{\mu ',f,\nu} \F[2]{a}{b}{g}{k}{\alpha,e,\gamma_{0}}{\gamma_{1},j,\eta} \F[2]{a}{j}{h}{d}{\eta,k,\epsilon_{0}}{\epsilon_{1},f,\nu}.
\end{align}
Among others, the PEPS and MPO tensors satisfy three consistency conditions derived from the left/right invertibility of $\FName{3}$ and the pentagon equations involving $\FName{2}$ and $\FName{3}$ respectively:
\begin{align}
	\begin{tikzarray}[scale=.15, every node/.style={scale=0.8}]{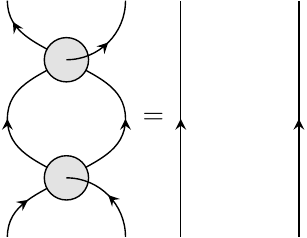}
		\begin{scope}[local bounding box=sc1]
			\path[] (-4.5,-8) rectangle (4.5,8);
			\draw [draw=DColor,mid arrow reversed] (0,4) to [out=-30,in=90] (4,0) to [out=270,in=30] (0,-4);
			\draw [draw=DColor,mid arrow reversed] (0,4) to [out=210,in=90] (-4,0) to [out=270,in=150] (0,-4);
			\draw [draw=DColor,TQ arrow] (0,4) to [out=150,in=270] (-4,8);
			\draw [draw=DColor,mid arrow reversed] (0,-4) to [out=-150,in=-270] (-4,-8);
			\draw [fill=gray!22] (0,4) circle (1.5);
			\draw [fill=gray!22] (0,-4) circle (1.5);
			\draw [draw=DColor,mid arrow] (0,4) to [out=0,in=-90] (4,8);
			\draw [draw=DColor,mid arrow reversed] (0,-4) to [out=0,in=90] (4,-8);
		\end{scope}
		\begin{scope}[local bounding box=sc2,anchor=west,shift={($(sc1.east)+(0,0)$)}]
			\node at (0,0) {$=$};
		\end{scope}
		\begin{scope}[local bounding box=sc3,anchor=west,shift={($(sc2.east)+(4.5,0)$)}]
			\path[] (-4.5,-8) rectangle (4.5,8);
			\draw [draw=DColor,mid arrow] (-4,-8) -- (-4,8);
			\draw [draw=DColor,mid arrow] (4,-8) -- (4,8);
		\end{scope}
	\end{tikzarray},\,
	\begin{tikzarray}[scale=.15, every node/.style={scale=0.8}]{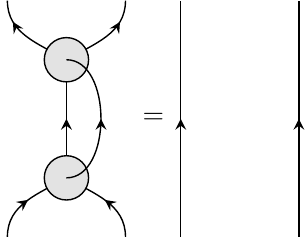}
		\begin{scope}[local bounding box=sc1]
			\path[] (-4.5,-8) rectangle (4.5,8);
			\draw [draw=DColor,TQ arrow] (0,4) to [out=30,in=270] (4,8);
			\draw [draw=DColor,TQ arrow] (0,4) to [out=150,in=270] (-4,8);
			\draw [draw=DColor,mid arrow reversed] (0,-4) to [out=-30,in=-270] (4,-8);
			\draw [draw=DColor,mid arrow reversed] (0,-4) to [out=-150,in=-270] (-4,-8);
			\draw [draw=DColor,mid arrow] (0,-4) -- (0,4);
			\draw[fill=gray!22] (0,4) circle (1.5);
			\draw[fill=gray!22] (0,-4) circle (1.5);
			\draw [draw=DColor,mid arrow reversed] (0,4) to [out=0,in=0] (0,-4);
		\end{scope}
		\begin{scope}[local bounding box=sc2,anchor=west,shift={($(sc1.east)+(0,0)$)}]
			\node at (0,0) {$=$};
		\end{scope}
		\begin{scope}[local bounding box=sc3,anchor=west,shift={($(sc2.east)+(4.5,0)$)}]
			\path[] (-4.5,-8) rectangle (4.5,8);
			\draw [draw=DColor,mid arrow] (-4,-8) -- (-4,8);
			\draw [draw=DColor,mid arrow] (4,-8) -- (4,8);
		\end{scope}
	\end{tikzarray},\,
	\begin{tikzarray}[scale=.15, every node/.style={scale=0.8}]{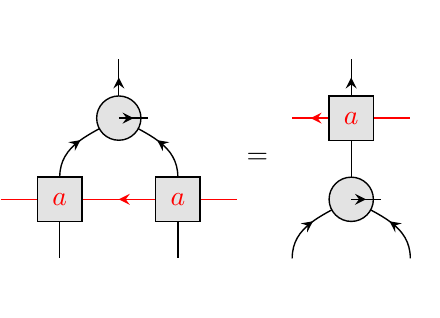}
		\begin{scope}[local bounding box=sc1]
			\path[] (-6,-8) rectangle (6,8);
			\draw [draw=DColor,mid arrow reversed] (0,0) to [out=-30,in=-270] (4,-4);
			\draw [draw=DColor,mid arrow reversed] (0,0) to [out=-150,in=-270] (-4,-4);
			\draw [draw=DColor,mid arrow] (0,1.5) -- (0,4);
			\draw [draw=CColor,mid arrow reversed] (-8,-5.5) -- (8,-5.5);
			\draw [draw=DColor] (-4,-7) -- (-4,-9.5);
			\draw [draw=DColor] (4,-7) -- (4,-9.5);
			\draw [fill=gray!22] (0,0) circle (1.5);
			\draw [draw=DColor,mid arrow] (0,0) -- (2,0);
			\draw [fill=gray!22] (-5.5,-7) rectangle (-2.5,-4);
			\draw [fill=gray!22] (2.5,-7) rectangle (5.5,-4);
			\node at (-4,-5.5) {$\catOb[C]{a}$};
			\node at (4,-5.5) {$\catOb[C]{a}$};
		\end{scope}
		\begin{scope}[local bounding box=sc2,anchor=west,shift={($(sc1.east)+(0,-2)$)}]
			\node at (0,0) {$=$};
		\end{scope}
		\begin{scope}[local bounding box=sc3,anchor=west,shift={($(sc2.east)+(5,-2.75)$)}]
			\path[] (-6,-8) rectangle (6,8);
			\draw [draw=DColor,mid arrow reversed] (0,0) to [out=-30,in=-270] (4,-4);
			\draw [draw=DColor,mid arrow reversed] (0,0) to [out=-150,in=-270] (-4,-4);
			\draw [draw=DColor] (0,1.5) -- (0,4);
			\draw [draw=DColor,mid arrow] (0,7) -- (0,9.5);
			\draw [draw=CColor,mid arrow reversed] (-4,5.5) -- (-1.5,5.5);
			\draw [draw=CColor] (1.5,5.5) -- (4,5.5);
			\draw [fill=gray!22] (0,0) circle (1.5);
			\draw [draw=DColor,mid arrow] (0,0) -- (2,0);
			\draw [fill=gray!22] (-1.5,4) rectangle (1.5,7);
			\node[anchor=center] at (0,5.5) {$\catOb[C]{a}$};
		\end{scope}
	\end{tikzarray}.
\end{align}
Using these, we can rewrite the MPO injectivity condition as
\begin{subalign}
	\frac{\dimM{e}\dimM{f}}{\dimD{c'}\FPdim{\D}}
	\begin{tikzarray}[scale=.15, every node/.style={scale=0.8}]{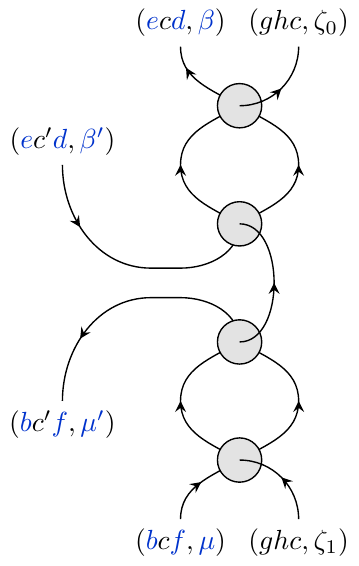}
		\draw [draw=DColor,mid arrow] (0,4) to [out=30,in=270] (4,8) to [out=90,in=-30] (0,12);
		\draw [draw=DColor,mid arrow] (0,4) to [out=150,in=270] (-4,8) to [out=90,in=210] (0,12);
		\draw [draw=DColor,TQ arrow] (0,12) to [out=150,in=270] (-4,16);
		\draw [draw=DColor,TQ arrow reversed] (0,4) to [out=270,in=0] (-4,1) -- (-6,1) to [out=180,in=270] (-12,8);
		\draw [draw=DColor,mid arrow reversed] (0,-4) to [out=-30,in=90] (4,-8) to [out=-90,in=30] (0,-12);
		\draw [draw=DColor,mid arrow reversed] (0,-4) to [out=-150,in=90] (-4,-8) to [out=-90,in=150] (0,-12);
		\draw [draw=DColor,mid arrow reversed] (0,-12) to [out=-150,in=-270] (-4,-16);
		\draw [draw=DColor,TQ arrow] (0,-4) to [out=90,in=0] (-4,-1) -- (-6,-1) to [out=180,in=90] (-12,-8);
		\draw[fill=gray!22] (0,4) circle (1.5);
		\draw[fill=gray!22] (0,-4) circle (1.5);
		\draw[fill=gray!22] (0,12) circle (1.5);
		\draw[fill=gray!22] (0,-12) circle (1.5);
		\draw [draw=DColor,mid arrow reversed] (0,4) to [out=0,in=0] (0,-4);
		\draw [draw=DColor,mid arrow] (0,12) to [out=0,in=-90] (4,16);
		\draw [draw=DColor,mid arrow reversed] (0,-12) to [out=0,in=90] (4,-16);
		\node[above] at (-12,8) {$(\catOb[M]{e}\catOb[D]{c'}\catOb[M]{d},\catOb[M]{\beta'})$};
		\node[above] at (-4,16) {$(\catOb[M]{e}\catOb[D]{c}\catOb[M]{d},\catOb[M]{\beta})$};
		\node[above] at (4,16) {$(\catOb[D]{g}\catOb[D]{h}\catOb[D]{c},\catOb[D]{\zeta_0})$};
		\node[below] at (-12,-8) {$(\catOb[M]{b}\catOb[D]{c'}\catOb[M]{f},\catOb[M]{\mu'})$};
		\node[below] at (-4,-16) {$(\catOb[M]{b}\catOb[D]{c}\catOb[M]{f},\catOb[M]{\mu})$};
		\node[below] at (4,-16) {$(\catOb[D]{g}\catOb[D]{h}\catOb[D]{c},\catOb[D]{\zeta_1})$};
	\end{tikzarray} & = \sum_{\catOb[C]{a}} \frac{\dimC{a}}{\FPdim{\C}}
	\begin{tikzarray}[scale=.15, every node/.style={scale=0.8}]{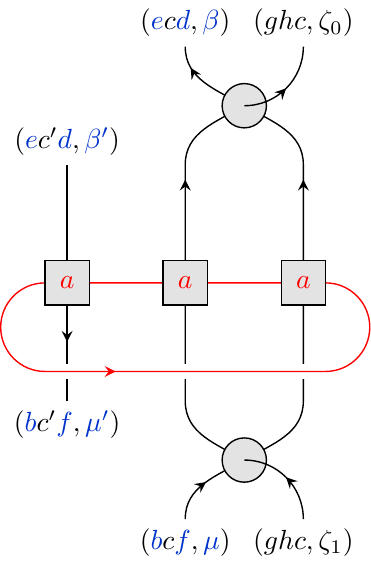}
		\draw [DColor,TQ arrow] (4,-12) to [out=150,in=-90] (0,-8) -- (0,8) to [out=90,in=210] (4,12);
		\draw [DColor,Q arrow reversed] (-8,-8) -- (-8,8);
		\draw [DColor,TQ arrow] (4,-12) to [out=30,in=-90] (8,-8) -- (8,8) to [out = 90,in=-30] (4,12);
		\draw [draw=DColor,TQ arrow] (4,12) to [out=150,in=270] (0,16);
		\draw [draw=DColor,mid arrow reversed] (4,-12) to [out=-150,in=-270] (0,-16);
		\draw[fill=gray!22] (4,12) circle (1.5);
		\draw[fill=gray!22] (4,-12) circle (1.5);
		\draw [draw=DColor,mid arrow] (4,12) to [out=0,in=-90] (8,16);
		\draw [draw=DColor,mid arrow reversed] (4,-12) to [out=0,in=90] (8,-16);
		\fill[white] (-0.5,-6.5) rectangle (0.5,-5.5);
		\fill[white] (-8.5,-6.5) rectangle (-7.5,-5.5);
		\fill[white] (7.5,-6.5) rectangle (8.5,-5.5);
		\draw[CColor,TQ arrow reversed] (-9.5,0) -- (9.5,0) arc (90:-90:3) -- (-9.5,-6) arc (270:90:3);
		\draw[fill=gray!22] (-1.5,-1.5) rectangle (1.5,1.5);
		\draw[fill=gray!22] (-9.5,-1.5) rectangle (-6.5,1.5);
		\draw[fill=gray!22] (6.5,-1.5) rectangle (9.5,1.5);
		\node at (0,0) {$\catOb[C]{a}$};
		\node at (8,0) {$\catOb[C]{a}$};
		\node at (-8,0) {$\catOb[C]{a}$};
		\node[above] at (-8,8) {$(\catOb[M]{e}\catOb[D]{c'}\catOb[M]{d},\catOb[M]{\beta'})$};
		\node[above] at (0,16) {$(\catOb[M]{e}\catOb[D]{c}\catOb[M]{d},\catOb[M]{\beta})$};
		\node[above] at (8,16) {$(\catOb[D]{g}\catOb[D]{h}\catOb[D]{c},\catOb[D]{\zeta_0})$};
		\node[below] at (-8,-8) {$(\catOb[M]{b}\catOb[D]{c'}\catOb[M]{f},\catOb[M]{\mu'})$};
		\node[below] at (0,-16) {$(\catOb[M]{b}\catOb[D]{c}\catOb[M]{f},\catOb[M]{\mu})$};
		\node[below] at (8,-16) {$(\catOb[D]{g}\catOb[D]{h}\catOb[D]{c},\catOb[D]{\zeta_1})$};
	\end{tikzarray}  \\
	\iff \frac{\dimM{e}\dimM{f}}{\dimD{c'}\FPdim{\D}}
	\begin{tikzarray}[scale=.15, every node/.style={scale=0.8}]{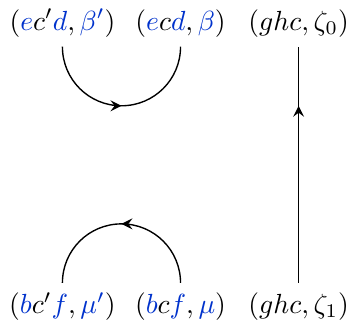}
		\path[] (-8,-8) rectangle (6.5,8);
		\draw [draw=DColor,mid arrow] (-8,8) to [out=-90,in=180] (-4,4) to [out=0,in=-90] (0,8);
		\draw [draw=DColor,mid arrow reversed] (-8,-8) to [out=90,in=180] (-4,-4) to [out=0,in=90] (0,-8);
		\draw [DColor,TQ arrow] (8,-8) -- (8,8);
		\node[above] at (-8,8) {$(\catOb[M]{e}\catOb[D]{c'}\catOb[M]{d},\catOb[M]{\beta'})$};
		\node[above] at (0,8) {$(\catOb[M]{e}\catOb[D]{c}\catOb[M]{d},\catOb[M]{\beta})$};
		\node[above] at (8,8) {$(\catOb[D]{g}\catOb[D]{h}\catOb[D]{c},\catOb[D]{\zeta_0})$};
		\node[below] at (-8,-8) {$(\catOb[M]{b}\catOb[D]{c'}\catOb[M]{f},\catOb[M]{\mu'})$};
		\node[below] at (0,-8) {$(\catOb[M]{b}\catOb[D]{c}\catOb[M]{f},\catOb[M]{\mu})$};
		\node[below] at (8,-8) {$(\catOb[D]{g}\catOb[D]{h}\catOb[D]{c},\catOb[D]{\zeta_1})$};
	\end{tikzarray}  & = \sum_{\catOb[C]{a}} \frac{\dimC{a}}{\FPdim{\C}}
	\begin{tikzarray}[scale=.15, every node/.style={scale=0.8}]{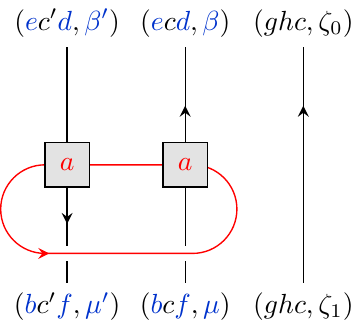}
		\path[] (-8,-8) rectangle (6.5,8);
		\draw [DColor,TQ arrow] (0,-8) -- (0,8);
		\draw [DColor,Q arrow reversed] (-8,-8) -- (-8,8);
		\draw [DColor,TQ arrow] (8,-8) -- (8,8);
		\fill[white] (-0.5,-6.5) rectangle (0.5,-5.5);
		\fill[white] (-8.5,-6.5) rectangle (-7.5,-5.5);
		\draw[CColor,TQ arrow reversed] (-9.5,0) -- (0.5,0) arc (90:-90:3) -- (-9.5,-6) arc (270:90:3);
		\draw[fill=gray!22] (-1.5,-1.5) rectangle (1.5,1.5);
		\draw[fill=gray!22] (-9.5,-1.5) rectangle (-6.5,1.5);
		\node at (0,0) {$\catOb[C]{a}$};
		\node at (-8,0) {$\catOb[C]{a}$};
		\node[above] at (-8,8) {$(\catOb[M]{e}\catOb[D]{c'}\catOb[M]{d},\catOb[M]{\beta'})$};
		\node[above] at (0,8) {$(\catOb[M]{e}\catOb[D]{c}\catOb[M]{d},\catOb[M]{\beta})$};
		\node[above] at (8,8) {$(\catOb[D]{g}\catOb[D]{h}\catOb[D]{c},\catOb[D]{\zeta_0})$};
		\node[below] at (-8,-8) {$(\catOb[M]{b}\catOb[D]{c'}\catOb[M]{f},\catOb[M]{\mu'})$};
		\node[below] at (0,-8) {$(\catOb[M]{b}\catOb[D]{c}\catOb[M]{f},\catOb[M]{\mu})$};
		\node[below] at (8,-8) {$(\catOb[D]{g}\catOb[D]{h}\catOb[D]{c},\catOb[D]{\zeta_1})$};
	\end{tikzarray},
\end{subalign}
the nontrivial part of which evaluates to
\begin{align}
	\deltaD{c}{c^\prime}\deltaM{\beta}{\beta^\prime}\deltaM{\mu}{\mu^\prime}  \frac{\dimM{e} \dimM{f}}{\dimD{c}} = \sum_{\substack{
	\catOb[C]{a} \\
			\catOb[M]{\alpha},\catOb[M]{\nu}}} \dimC{a} \F[2]{a}{b}{c}{d}{\alpha,e,\beta}{\mu,f,\nu} \Fi[2]{a}{b}{c'}{d}{\alpha,e,\beta '}{\mu ',f,\nu}.
\end{align}
This is exactly the Schur orthogonality condition \cref{eqn:matrixElementOrthog}. Together with the dimensionality condition, it guarantees that the fusion categories $\C$ and $\D$ are Morita equivalent (i.e., $\bimodule{\C}{\M}{\D}$ is invertible), or equivalently, they have tensor-equivalent Drinfel'd centers $\drinfeld{\C} \tensorEquiv \drinfeld{\D}$~\cite{Mueger2003,Mueger2003a,Etingof2010,Etingof2011}.

\section{Remarks}\label{sec:remarks}

The Schur orthogonality conditions are a cornerstone in the representation theory of finite groups and as such play an important role in physics. A particular incarnation is known as the Wigner-Eckart theorem, stating that a generic $G$-symmetric tensor operator can be decomposed into Clebsch-Gordan coefficients. Recently there has been much interest in a generalized notion of symmetries, which in (1+1)-dimensions is described by fusion categories. In lattice models, these symmetries are represented as MPOs, which in turn are described by bimodule categories. In this context, the generalized Schur orthogonality conditions play the same role as for ordinary finite group symmetries, and they allow one to show a generalized Wigner-Eckart theorem, the proof of which parallels the derivation of the MPO-injectivity condition. The generalized Wigner-Eckart theorem allows one to interpret the tensors \cref{eq:PEPSdef} as generalized Clebsch-Gordan coefficients, an insight which has recently been exploited to provide a comprehensive picture of dualities in (1+1)-dimensions~\cite{Lootens2022}.

Generalizations of Schur orthogonality conditions from finite groups to compact groups proceed in a fairly straightforward way by replacing the sum over group elements to a Haar integral over the group, and follow from the Peter-Weyl theorem. Many physical systems have symmetries described by compact Lie groups, and it would be interesting to investigate the orthogonality conditions for invertible bimodule categories that go beyond the standard representation theory case, discussed in \cref{app:VecG} for finite groups.

A further generalization is obtained by categorification, in which case the central question we addressed in this work amounts to asking what the conditions are on the skeletal data of a bimodule 2-category for it to be invertible. There are several difficulties associated to this, the most obvious being that the skeletal approach for fusion 2-categories is not yet well established. Furthermore, to generalize our proof of the orthogonality conditions one would require the notion of a weak Hopf 2-algebra and its representation theory, which again has not yet been worked out in the literature.
\acknowledgments

We thank Clement Delcamp, J\"urgen Fuchs, Peter Huston, Corey Jones, Andr\'as Moln\'ar, Christoph Schweigert and Tim Seynnaeve for useful discussions.
L.L.~is supported by a PhD fellowship from the Research Foundation Flanders (FWO).

This work was initiated at the workshop “\href{https://aimath.org/pastworkshops/fusiontensorV.html}{Fusion categories and tensor networks}”, hosted by the American Institute
of Mathematics in March 2021. We thank AIM for their generosity.

\bibliographystyle{apsrev_jacob}
\bibliography{refs}

\appendix


\section{Notation}\label{sec:symbols}

\renewcommand{\arraystretch}{1.1}
\setlength{\tabcolsep}{20pt}

\begin{tabular}{!{\color[gray]{.8}\vrule}>{}c!{\color[gray]{.8}\vrule}l!{\color[gray]{.8}\vrule}}
	\greyhline
	$\irr{\C}$               & The set of simple objects of the category $\C$                                                                    \\
	$\rk{\C}$                & The number of (isomorphism classes of) simple objects of the category $\C$                                        \\
	$\dual{a}$               & The object dual to $a\in \irr{\C}$                                                                                \\
	$\homSpace{\C}{a}{b}$    & The space of morphisms from $a$ to $b$ in the category $\C$                                                       \\
	$\sum_{a}$               & Shorthand for $\sum_{a\in\irr{\cat{X}}}$, where $\cat{X}$ is the appropriate category                             \\
	$\N{a}{b}{c}$            & Vector space dimension of $\homSpace{\C}{a\otimes b}{c}$.                                                         \\
	$e_{ij}^{\alpha}$        & Matrix unit, multiplying as $e_{ij}^{\alpha}e_{kl}^{\beta}=\delta_{j}^{k}\delta_{\alpha}^{\beta} e_{il}^{\alpha}$ \\
	\greyhline
	\multicolumn{2}{c}{Dimensions}                                                                                                               \\
	\greyhline
	$\dimQ{X}$               & Frobenius-Perron dimension of object $X$                                                                          \\
	$\FPdim{\C}$             & Frobenius-Perron dimension of category $\C$                                                                       \\
	$\dim{V}$                & Dimension of $V$ as a $\mathbb{C}$-vector space                                                                   \\
	\greyhline
	\multicolumn{2}{c}{Associators}                                                                                                              \\
	\greyhline
	$\FName{0}$              & $F$-symbol in category $\C$                                                                                       \\
	$\FName{1}$              & Left module associator of $\leftmodule{\C}{\M}$                                                                   \\
	$\FName{2}$              & Bimodule associator of $\bimodule{\C}{\M}{\D}$                                                                    \\
	$\FName{3}$              & Right module associator of $\rightmodule{\M}{\D}$                                                                 \\
	$\FName{4}$              & $F$-symbol in category $\D$                                                                                       \\
	\greyhline
	\multicolumn{2}{c}{Weak Hopf algebra (See \onlinecite{Boehm1999} for definitions)}                                                           \\
	\greyhline
	$\dualSpace{A}$          & Dual space $\Hom(A,\mathbb{C})$                                                                                   \\
	$\mu$                    & Product (usually written as juxtaposition)                                                                        \\
	$\eta$                   & Unit                                                                                                              \\
	$\Delta$                 & Coproduct                                                                                                         \\
	$\varepsilon$            & Counit                                                                                                            \\
	$S$                      & Antipode                                                                                                          \\
	$\Lambda$                & Haar integral                                                                                                     \\
	$\lambda$                & Haar measure in $\dualSpace{A}$                                                                                   \\
	$x_{(1)}\otimes x_{(2)}$ & Sweedler notation for $\Delta x$                                                                                  \\
	\greyhline
	\multicolumn{2}{c}{Representations and modules}                                                                                              \\
	\greyhline
	$a\cdot v$               & Action of $a$ on $v$                                                                                              \\
	$\rho_x(a)$              & Map in $\mathrm{End}{(V_x)}$ which implements $a\cdot v$                                                          \\
	$\chi_x$                 & Character in representation labeled by $x$                                                                        \\
	$V\boxtimes W$           & Product on modules $\set{x\in V\otimes W}{\Delta(1)x=x}$                                                          \\
	\greyhline
\end{tabular}

\section{The weak Hopf algebra structure of \texorpdfstring{$\ann{\C}{\M}$}{Ann(M)}}\label{app:WHA}

Given a unitary module category $\leftmodule{\C}{\M}$, the module category $\ann{\C}{\M}$ can be equipped with the structure of a weak Hopf algebra~\cite{Kitaev2012}. We refer to \onlinecite{Boehm1999} for a definition of weak Hopf algebras.

In this section, we review this, providing explicit actions of all maps in terms of the `picture basis'. The underlying vector space is
\begin{align}
	A & =\Cspan
	\set{
		\tub{a,b}{c,d}{\alpha,x,\beta}
	}{
		\catOb[M]{a},\catOb[M]{b},\catOb[M]{c},\catOb[M]{d}\in\irr{\M},\catOb[C]{x}\in\irr{\C},1\leq\textcolor{MColor}{\alpha}\leq \NCM{x}{a}{c},1\leq\textcolor{MColor}{\beta}\leq \NCM{x}{b}{d}
	},
\end{align}
which has dimension
\begin{align}
	\dim(A)=\sum_{\substack{\catOb[M]{a},\catOb[M]{b},\catOb[M]{c},\catOb[M]{d}\in\irr{\M} \\\catOb[C]{x}\in\irr{\C}}} \NCM{x}{a}{c}\NCM{x}{b}{d}.
\end{align}

Multiplication is defined by concentrically stacking the diagrams (introduced in \cref{eqn:pictureBasis}). Such a product can only be nonzero if the marked boundary points match. The associator isomorphisms can be used to reduce the composite diagram to the picture basis:
\begin{align}
	\mu:\tub{a^\prime,b^\prime}{c^\prime,d^\prime}{\alpha^\prime,x^\prime,\beta^\prime}\otimes\tub{a,b}{c,d}{\alpha,x,\beta}
	\mapsto
	\delta_{\catOb[M]{c}}^{\catOb[M]{a^\prime}}\delta_{\catOb[M]{d}}^{\catOb[M]{b^\prime}}
	\sum_{\substack{\catOb[0]{y},\catOb[0]{\zeta} \\\catOb[1]{\mu},\catOb[1]{\nu}}}\sqrt{\frac{\dimC{x}\dimC{x^\prime}}{\dimC{y}}}
	\F[1]{x^\prime}{x}{a}{c^\prime}{\zeta,y,\mu}{\alpha,c,\alpha^\prime} \Fi[1]{x^\prime}{x}{b}{d^\prime}{\zeta,y,\nu}{\beta,d,\beta^\prime}
	\tub{a,b}{c^\prime,d^\prime}{\mu,y,\nu}.
\end{align}
The unit on $A$ is given by
\begin{align}
	\eta: 1\mapsto\one = \sum_{\catOb[M]{a},\catOb[M]{b}}  \tub{a,b}{a,b}{1,1,1}.
\end{align}
With these maps, $(A,\mu,\eta)$ is an associative, unital algebra.

Additionally, we can equip $\ann{\C}{\M}$ with a coproduct
\begin{align}
	\Delta: \tub{a,b}{c,d}{\alpha,x,\beta}\mapsto \frac{1}{\sqrt{\dimC{x}}}\sum_{\catOb[M]{e},\catOb[M]{f},\catOb[M]{\mu}}\tub{e,b}{f,d}{\mu,x,\beta}\otimes \tub{a,e}{c,f}{\alpha,x,\mu},
\end{align}
and counit
\begin{align}
	\varepsilon:\tub{a,b}{c,d}{\alpha,x,\beta}\mapsto \delta_{\catOb[1]{\alpha}}^{\catOb[1]{\beta}} \delta_{\catOb[1]{a}}^{\catOb[1]{b}} \delta_{\catOb[1]{c}}^{\catOb[1]{d}} \sqrt{\dimC{x}}.
\end{align}
The triple $(A,\Delta,\varepsilon)$ is a coassociative, counital algebra, and $(A,\mu,\eta,\Delta,\varepsilon)$ is a weak bialgebra.

The antipode on $\ann{\C}{\M}$ is defined by
\begin{align}
	S:\tub{a,b}{c,d}{\alpha,x,\beta}\mapsto \frac{\dimM{b} \dimC{x}}{\dimM{d}} \sum _{\catOb[1]{\mu},\catOb[1]{\nu}}\F[1]{\dual{x}}{x}{a}{a}{1,1,1}{\alpha,c,\nu} \Fi[1]{\dual{x}}{x}{b}{b}{1,1,1}{\beta,d,\mu} \tub{d,c}{b,a}{\mu,\dual{x},\nu}.
\end{align}
We remark that $\ann{\C}{\M}$ becomes a true Hopf algebra, exactly when $\M$ has a single simple object.

Since the underlying module category is unitary, $\ann{\C}{\M}$ is a $*$-WHA, with
\begin{align}
	*:\tub{a,b}{c,d}{\alpha,x,\beta}\mapsto \dimC{x}\sqrt{\frac{\dimM{a}\dimM{b}}{\dimM{c}\dimM{d}}} \sum _{\catOb[1]{\mu},\catOb[1]{\nu}} \F[1]{\dual{x}}{x}{b}{b}{1,1,1}{\beta,d,\nu} \Fi[1]{\dual{x}}{x}{a}{a}{1,1,1}{\alpha,c,\mu} \tub{c,d}{a,b}{\mu,\dual{x},\nu}.
\end{align}

By manipulating string diagrams, it can be readily verified that
\begin{subalign}
	\Lambda:=\sum_{\substack{\catOb[M]{a},\catOb[M]{b} \\\catOb[C]{x},\catOb[M]{\alpha}}} \frac{\sqrt{\dimC{x}}}{\dimM{a} \dimM{b} \rk{\M}}  \tub{a,a}{b,b}{\alpha,x,\alpha}\\
	\lambda:\tub{a,b}{c,d}{\alpha,x,\beta}\mapsto \delta_{\catOb[M]{\alpha}}^{1} \delta_{\catOb[M]{\beta}}^{1} \delta_{\catOb[C]{x}}^{1} \delta_{\catOb[M]{a}}^{\catOb[M]{c}} \delta_{\catOb[M]{b}}^{\catOb[M]{d}}  \rk{\M} \,\dimM{a}^2,
\end{subalign}
are a dual pair of Haar integrals. By \onlinecite[Theorem 3.27]{Boehm1999}, this establishes semisimplicity of $\ann{\C}{\M}$. Finally, $\ann{\C}{\M}$ has a faithful $*$-representation, namely the regular representation equipped with the inner product $\langle x,y\rangle_A=\lambda(x^*y)$, as so is a $C^*$-WHA.

In addition to the Haar integrals, we need the \define{canonical grouplike element}~\cite{Boehm1999}, and its inverse
\begin{subalign}[eqn:grouplike]
	g=       & \sum_{\catOb[M]{a},\catOb[M]{b}} \frac{\dimM{a}}{\dimM{b}} \tub{a,b}{a,b}{1,1,1}, \\
	g^{-1} = & \sum_{\catOb[M]{a},\catOb[M]{b}} \frac{\dimM{b}}{\dimM{a}} \tub{a,b}{a,b}{1,1,1}.
\end{subalign}


\section{Example: \texorpdfstring{$\vvec{G}$}{VecG}}\label{app:VecG}

Let $G$ be a finite group. The category of $G$-graded (finite dimensional) vector spaces is denoted $\vvec{G}$. The simple objects are one-dimensional vector spaces in degree $g$, for each $g\in G$. By a slight abuse of notation, we will denote these objects by their group label. Fusion rules for $\vvec{G}$ are given by the group multiplication, so the valid trivalent vertices for simple objects are
\begin{align}
	\begin{tikzarray}[scale=0.4, every node/.style={scale=0.8}]{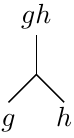}
		\trivalentvertex{gh}{g}{h}{};
	\end{tikzarray},
\end{align}
all with multiplicity $1$. All $\FName{0}$-symbols are $+1$ when allowed by the fusion rules.

For all $G$, the category of finite dimensional vector spaces (denoted $\vvec{}$) can be equipped with the structure of a module category over $\vvec{G}$. We will denote the single simple object of $\vvec{}$ by $\catOb[1]{\vecOb}$. All fusion spaces $\vvec{}\left(g\triangleright\catOb[1]{\vecOb},\catOb[1]{\vecOb}\right)$ are one-dimensional, and all $\FName{1}$-symbols are $+1$. The module object has Frobenius-Perron dimension $\dimM{\vecOb}= \sqrt{|G|}$.

For groups with non-trivial projective representations, we can twist the group action, giving $\FName{1}$ symbols
\begin{align}
	\F[1]{g}{h}{\vecOb}{\vecOb}{k}{\vecOb} & = \delta_{gh}^{k}\, \phi(g,h),\tag{Twisted group action}
\end{align}
where $\phi\in H^{2}(G,\mathbb{C}^{\times})$ is a normalized 2-cocycle.

With this data fixed, the first pentagon equation governing $\FName{2}$ is given by
\begin{align}
	\sum_{\zeta } \F[2]{g}{\vecOb}{\pi}{\vecOb}{1,\vecOb,\alpha}{\zeta,\vecOb,1} \F[2]{h}{\vecOb}{\pi}{\vecOb}{1,\vecOb,\zeta}{\beta,\vecOb,1}=\F[2]{g h}{\vecOb}{\pi}{\vecOb}{1,\vecOb,\alpha}{\beta,\vecOb,1},
\end{align}
which specifies that $\F[2]{g}{\vecOb}{\pi}{\vecOb}{1,\vecOb,\alpha}{\zeta,\vecOb,1}$ are the matrix elements of a representation $\pi$ of $G$
\begin{align}
	\F[2]{g}{\vecOb}{\pi}{\vecOb}{1,\vecOb,\alpha}{\beta,\vecOb,1} & =\bigl[\rho_{\pi}(g)\bigr]_{\alpha\beta}.
\end{align}
For this special case, our results reduce to familiar expressions. Invertibility of $\vvec{}$ (\cref{eqn:invertibilityConditions}) occurs when
\begin{align}
	\frac{1}{|G|} \sum_{g\in G}\chi_{c}(g)\overline{\chi_{c'}(g)} & \oldequiv \deltaD{c}{c'},
\end{align}
which is Schur's orthogonality relation for (irreducible) characters of finite groups. In particular, this shows that the dual category $(\vvec{G})^{*}_{\vvec{}}$ is $\rrep{G}$ as expected.

The matrix element condition of \cref{thm:orthogMatrixElements} also reduces to the familiar expression
\begin{align}
	\sum_{g\in G}\bigl[\rho_{c}(g)\bigr]_{\alpha\beta}\overline{\bigl[\rho_{c}(g)\bigr]}_{\alpha'\beta'} & = \frac{|G|}{d_c}\deltaD{c}{c'}\delta_{\alpha}^{\alpha'}\delta_{\beta}^{\beta'} = \frac{|G|}{\dim{V_c}}\deltaD{c}{c'}\delta_{\alpha}^{\alpha'}\delta_{\beta}^{\beta'},
\end{align}
for orthogonality of matrix elements of irreducible representations.

The intertwining maps expressing the tensor product of representations in terms of irreps solve the next pentagon equation
\begin{align}
	\bigl[\rho_{\sigma}\left(g\right)\bigr]_{\alpha\alpha '} \bigl[\rho_{\tau}\left(g\right)\bigr]_{\beta\beta '} & =\sum _{\substack{\pi\in\sigma\otimes\tau \\\mu,\nu,\epsilon }} \F[3]{\vecOb}{\sigma}{\tau}{\vecOb}{\alpha,\vecOb,\beta}{\mu,\pi,\nu} \bigl[\rho_{\pi}\left(g\right)\bigr]_{\nu\epsilon} \Fi[3]{\vecOb}{\sigma}{\tau}{\vecOb}{\alpha ',\vecOb,\beta '}{\mu,\pi,\epsilon},\label{eqn:CGs}
\end{align}
which recover the Clebsch-Gordan coefficients
\begin{align}
	\F[3]{\vecOb}{a}{b}{\vecOb}{\alpha,\vecOb,\beta}{\mu,c,\nu} & =\CG[\mu]{a}{b}{c}{\alpha}{\beta}{\nu} \tag{Clebsch-Gordan coefficients}.
\end{align}
Finally, the remaining skeletal data is
\begin{align}
	\F[4]{a}{b}{c}{d}{\alpha,e,\beta}{\mu,f,\nu} & =W. \tag{Racah W-coefficients/$6j$-symbols}
\end{align}

\section{Weak Hopf algebra}

A weak Hopf algebra is a vector space $A$ equipped with weakly compatible algebra and coalgebra structures. For convenience, we review these compatibility relations here.

\subsubsection{Algebra}

$(A,\mu:A\otimes A\to A,\eta:\mathbb{C}\to A)$

\[
	\begin{tikzcd}
		{A\otimes A\otimes A} & {A\otimes A} \\
		{A\otimes A} & A
		\arrow["\mu\otimes\id", from=1-1, to=1-2]
		\arrow["\id\otimes\mu"', from=1-1, to=2-1]
		\arrow["\mu"', from=2-1, to=2-2]
		\arrow["\mu", from=1-2, to=2-2]
	\end{tikzcd}
	\hspace{2cm}
	\begin{tikzcd}
		A & {A\otimes A} \\
		{A\otimes A} & A
		\arrow["\eta\otimes\id", from=1-1, to=1-2]
		\arrow["\id\otimes\eta"', from=1-1, to=2-1]
		\arrow["\mu"', from=2-1, to=2-2]
		\arrow["\mu", from=1-2, to=2-2]
		\arrow["\id", from=1-1, to=2-2]
	\end{tikzcd}
\]

\subsubsection{Coalgebra}

$(A,\Delta:A\to A\otimes A,\varepsilon:A\to\mathbb{C})$

\[\begin{tikzcd}
		A & {A\otimes A} \\
		{A\otimes A} & {A\otimes A\otimes A}
		\arrow["\Delta"', from=1-1, to=2-1]
		\arrow["\Delta", from=1-1, to=1-2]
		\arrow["\Delta\otimes\id", from=1-2, to=2-2]
		\arrow["\id\otimes\Delta"', from=2-1, to=2-2]
	\end{tikzcd}
	\hspace{2cm}
	\begin{tikzcd}
		A & {A\otimes A} \\
		{A\otimes A} & A
		\arrow["\id", from=1-1, to=2-2]
		\arrow["\Delta", from=1-1, to=1-2]
		\arrow["\Delta"', from=1-1, to=2-1]
		\arrow["\id\otimes\varepsilon"', from=2-1, to=2-2]
		\arrow["\varepsilon\otimes\id", from=1-2, to=2-2]
	\end{tikzcd}\]

\subsubsection{Weak bialgebra}

$(A,\mu,\eta,\Delta,\varepsilon)$

In the following, $\sigma:A\otimes A\to A\otimes A$ swaps the factors.

\[\begin{tikzcd}
		{A\otimes A} & A & {A\otimes A} \\
		\\
		{A\otimes A\otimes A\otimes A} && {A\otimes A\otimes A\otimes A}
		\arrow["\Delta\otimes\Delta"', from=1-1, to=3-1]
		\arrow["\mu", from=1-1, to=1-2]
		\arrow["\Delta", from=1-2, to=1-3]
		\arrow["\mu\otimes\mu"', from=3-3, to=1-3]
		\arrow["\id\otimes\sigma\otimes\id"', from=3-1, to=3-3]
	\end{tikzcd}
\]
\[\begin{tikzcd}
		& {A\otimes A\otimes A\otimes A} && {A\otimes A\otimes A\otimes A} && {A\otimes A} \\
		{A\otimes A\otimes A} && {A\otimes A} && A && {\mathbb{C}} \\
		& {A\otimes A\otimes A\otimes A} &&&& {A\otimes A}
		\arrow["\mu\otimes\id", from=2-1, to=2-3]
		\arrow["\mu", from=2-3, to=2-5]
		\arrow["\varepsilon", from=2-5, to=2-7]
		\arrow["\id\otimes\Delta\otimes\id", from=2-1, to=1-2]
		\arrow["\id\otimes\sigma\otimes\id", from=1-2, to=1-4]
		\arrow["\mu\otimes\mu", from=1-4, to=1-6]
		\arrow["\varepsilon\otimes\varepsilon", from=1-6, to=2-7]
		\arrow["\varepsilon\otimes\varepsilon"', from=3-6, to=2-7]
		\arrow["\id\otimes\Delta\otimes\id"', from=2-1, to=3-2]
		\arrow["\mu\otimes\mu"', from=3-2, to=3-6]
	\end{tikzcd}\]
\[\begin{tikzcd}
		& {A\otimes A\otimes A\otimes A} && {A\otimes A\otimes A\otimes A} && {A\otimes A} \\
		{A\otimes A\otimes A} && {A\otimes A} && A && {\mathbb{C}} \\
		& {A\otimes A\otimes A\otimes A} &&&& {A\otimes A}
		\arrow["\Delta\otimes\id"', from=2-3, to=2-1]
		\arrow["\Delta"', from=2-5, to=2-3]
		\arrow["\eta"', from=2-7, to=2-5]
		\arrow["\id\otimes\mu\otimes\id"', from=1-2, to=2-1]
		\arrow["\id\otimes\sigma\otimes\id"', from=1-4, to=1-2]
		\arrow["\Delta\otimes\Delta"', from=1-6, to=1-4]
		\arrow["\eta\otimes\eta"', from=2-7, to=1-6]
		\arrow["\eta\otimes\eta", from=2-7, to=3-6]
		\arrow["\id\otimes\mu\otimes\id", from=3-2, to=2-1]
		\arrow["\Delta\otimes\Delta", from=3-6, to=3-2]
	\end{tikzcd}\]

\subsubsection{Weak Hopf algebra}

$(A,\mu,\eta,\Delta,\varepsilon,S:A\to A)$

\[\begin{tikzcd}
		{A\otimes A} &&& {A\otimes A} \\
		A &&& A \\
		{A\otimes A} & {A\otimes A\otimes A} & {A\otimes A\otimes A} & {A\otimes A}
		\arrow["\Delta", from=2-1, to=1-1]
		\arrow["{\id\otimes S}", from=1-1, to=1-4]
		\arrow["\mu", from=1-4, to=2-4]
		\arrow["\eta\otimes\id"', from=2-1, to=3-1]
		\arrow["\Delta\otimes\id"', from=3-1, to=3-2]
		\arrow["\id\otimes\sigma"', from=3-2, to=3-3]
		\arrow["\mu\otimes\id"', from=3-3, to=3-4]
		\arrow["\varepsilon\otimes\id"', from=3-4, to=2-4]
		\arrow["{\sqcap^{L}}"{description}, dashed, from=2-1, to=2-4]
	\end{tikzcd}\]
\[\begin{tikzcd}
		{A\otimes A} && {} & {A\otimes A} \\
		A &&& A \\
		{A\otimes A} & {A\otimes A\otimes A} & {A\otimes A\otimes A} & {A\otimes A}
		\arrow["\Delta", from=2-1, to=1-1]
		\arrow["{S\otimes \id}", from=1-1, to=1-4]
		\arrow["\mu", from=1-4, to=2-4]
		\arrow["\id\otimes\eta"', from=2-1, to=3-1]
		\arrow["\id\otimes\Delta"', from=3-1, to=3-2]
		\arrow["\sigma\otimes\id"', from=3-2, to=3-3]
		\arrow["\id\otimes\mu"', from=3-3, to=3-4]
		\arrow["\id\otimes\varepsilon"', from=3-4, to=2-4]
		\arrow["{\sqcap^{R}}"{description}, dashed, from=2-1, to=2-4]
	\end{tikzcd}\]
\[\begin{tikzcd}
		&& {} \\
		& A & A \\
		{A\otimes A} & {A\otimes A\otimes A} & {A\otimes A\otimes A} & {A\otimes A}
		\arrow["\Delta"', from=2-2, to=3-1]
		\arrow["\Delta\otimes\id"', from=3-1, to=3-2]
		\arrow["{S\otimes\id\otimes S}"', from=3-2, to=3-3]
		\arrow["\mu\otimes\id"', from=3-3, to=3-4]
		\arrow["\mu"', from=3-4, to=2-3]
		\arrow["S", from=2-2, to=2-3]
	\end{tikzcd}\]

Hopf if $S(\eta(1))=\eta(1)\otimes \eta(1)$.

\subsubsection{Weak Hopf *-algebra}

$(A,\mu,\eta,\Delta,\varepsilon,S,*:A\to A)$

\[\begin{tikzcd}
		{A\otimes A} & {} & A \\
		{A\otimes A} & {A\otimes A} & A
		\arrow["\mu", from=1-1, to=1-3]
		\arrow["\sigma"', from=1-1, to=2-1]
		\arrow["{*\otimes *}"', from=2-1, to=2-2]
		\arrow["{*}", from=1-3, to=2-3]
		\arrow["\mu"', from=2-2, to=2-3]
	\end{tikzcd}
	\hspace{2cm}
	\begin{tikzcd}
		{} & A & A \\
		&& A
		\arrow["{*}", from=1-2, to=1-3]
		\arrow["{*}", from=1-3, to=2-3]
		\arrow["\id"', from=1-2, to=2-3]
	\end{tikzcd}
	\hspace{2cm}
	\begin{tikzcd}
		{} & A & A \\
		& {A\otimes A} & {A\otimes A}
		\arrow["{*}", from=1-2, to=1-3]
		\arrow["\Delta", from=1-3, to=2-3]
		\arrow["\Delta"', from=1-2, to=2-2]
		\arrow["{*\otimes*}"', from=2-2, to=2-3]
	\end{tikzcd}
\]

\end{document}